\documentclass[11pt, a4paper]{article}
\usepackage{amssymb}
\usepackage{amsmath}

\numberwithin{equation}{section}

\usepackage{graphicx,color}
\usepackage{wrapfig,framed}

\newcommand{\diag}{{\rm diag}\hskip 0.5truemm}

\newtheorem{theorem}{Theorem}

\newtheorem{prop}[theorem]{Proposition}

\newtheorem{lemma}[theorem]{Lemma}

\newtheorem{cor}[theorem]{Corollary}

\newtheorem{remark}[theorem]{Remark}
\newtheorem{defi}[theorem]{Definition}

\def \= {\;=\;}
\def \+ {\,+\,}

\newcommand{\eqa}{\begin{eqnarray}}
\newcommand{\eeqa}{\end{eqnarray}}
\newcommand{\beq}{\begin{equation}}
\newcommand{\eeq}{\end{equation}}
\newcommand{\nn}{\nonumber}

\newcommand{\pal}{\partial}

\newcommand{\CC}{\mathbb{C}}
\newcommand{\QQ}{\mathbb{Q}}
\newcommand{\ZZ}{\mathbb{Z}}
\newcommand{\CW}{{\mathcal W}}
\newcommand{\CY}{{\mathcal Y}}

\newcommand{\bpsi}{{\boldsymbol \psi}}
\newcommand{\bla}{{\boldsymbol \lambda}}
\newcommand{\bom}{{\boldsymbol \Omega}}

\newcommand{\w}{\omega}

\newcommand{\pf}{\noindent{\it Proof: \ }}
\newcommand{\tr}{{\rm tr}}
\def\res{\mathop{\rm res}\limits}

\newcommand{\bt}{{\bf t}}
\newcommand{\by}{{\bf y}}

\newcommand{\bv}{{\bf V}}
\newcommand{\bu}{{\bf u}}
\newcommand{\bx}{{\bf x}}
\newcommand{\ii}{\sqrt{-1}}

\newcommand{\epf}{$\quad$\hfill
\raisebox{0.11truecm}{\fbox{}}\par\vskip0.4truecm}

\usepackage{array}
\newcolumntype{M}[1]{>{\centering\arraybackslash}m{#1}}
\newcolumntype{R}[1]{>{\raggedleft\arraybackslash}m{#1}}
\newcolumntype{N}{@{}m{0pt}@{}}

\def\Xint#1{\mathchoice
{\XXint\displaystyle\textstyle{#1}}%
{\XXint\textstyle\scriptstyle{#1}}%
{\XXint\scriptstyle\scriptscriptstyle{#1}}%
{\XXint\scriptscriptstyle\scriptscriptstyle{#1}}%
\!\int}
\def\XXint#1#2#3{{\setbox0=\hbox{$#1{#2#3}{\int}$}
\vcenter{\hbox{$#2#3$}}\kern-.5\wd0}}

\def\dashint{\Xint-}

\begin{document}

\title{Algebraic spectral curves over $\mathbb Q$ and their tau-functions}\par
\author{Boris DUBROVIN}

\maketitle

\begin{abstract}  Let $W(z)$ be a $n\times n$ matrix polynomial with rational coefficients. Denote $C$ the {\it spectral curve} $\det \left( w\cdot{\bf 1}-W(z)\right) =0$. Under some natural assumptions about the structure of $W(z)$ we prove that certain combinations of logarithmic derivatives of the Riemann theta-function of $C$ of an arbitrary order starting from the third one all take rational values at the point of the Jacobi variety $J(C)$ specified by the line bundle of eigenvectors of $W(z)$.
\end{abstract}

\section{Introduction}

Let
\beq\label{mapol}
W(z)=B^0 z^m+B^1 z^{m-1}+\dots+B^m, \quad B^i\in Mat_n(\mathbb C), \quad B^0 =\diag(b^0_1,\dots,b^0_n)
\eeq
be a $n\times n$ matrix polynomial.
Consider an algebraic curve $C$ defined by the characteristic equation
\beq\label{kurva}
R(z,w):=\det\left( w\cdot {\bf 1}-W(z)\right)=0.
\eeq
It will be called \emph{algebraic spectral curve} associated with the matrix-valued polynomial $W(z)$. Consider the generic situation when  the leading coefficient $B
^0$ has pairwise distinct eigenvalues.  In that case the Riemann surface $C\xrightarrow{z}{\bf P}^1$ has $n$ distinct infinite points $P_1\cup\dots \cup P_n =z^{-1}(\infty)$ labelled by the eigenvalues. It will also be assumed that the affine part of the curve \eqref{kurva} is smooth irreducible.  In sequel only such spectral curves will be considered. The genus $g$ of such a spectral curve is uniquely determined by the numbers $m$ and $n$, see eq. \eqref{nwg} below.

Assuming $g>0$ choose a canonical basis of cycles $a_1$, \dots, $a_g$, $b_1$, \dots, $b_g\in H_1(C, \mathbb Z)$,
\beq\label{canon}
a_i\circ a_j=b_i\circ b_j=0, \quad a_i\circ b_j=\delta_{ij}, \quad i, \, j=1, \dots, g.
\eeq
Let $\w_1, \dots \w_g$
be the basis of holomorphic differentials on $C$ normalized by the conditions
\beq\label{holo}
\oint_{a_j} \omega_k =2\pi \sqrt{-1}\, \delta_{jk}, \quad j, \, k=1, \dots, g.
\eeq
Denote
\beq\label{matper}
B_{jk}=\oint_{b_j} \omega_k, \quad j, \, k=1, \dots, g
\eeq
and let
\beq\label{theta}
\theta(\bu)=\sum_{{\bf n}\in \mathbb Z^g} \exp\left\{ \frac12\langle {\bf n}, B{\bf n}\rangle +\langle {\bf n},\bu\rangle \right\}
\eeq
be the Riemann theta-function of the curve $C$ associated with the chosen basis of cycles. Here $\bu=(u_1, \dots, u_g)$ is the vector of independent complex variables $u_1$, \dots, $u_g$, ${\bf n}=(n_1, \dots, n_g)\in\mathbb Z^g$,
$$
\langle {\bf n}, B{\bf n} \rangle=\sum_{j, \, k=1}^g B_{jk}n_j n_k, \quad \langle {\bf n}, {\bf u}\rangle=\sum_{k=1}^g n_k u_k.
$$
Also recall that the Jacobi variety, or simply Jacobian of the curve $C$ is defined as the quotient
$$
J(C)=\CC^g/\{ 2\pi \ii\, {\bf m} + B\, {\bf n}\, |\, {\bf m}, \, {\bf n}\in\mathbb Z^g\}.
$$

There is a natural line bundle of degree $g+n-1$ on the spectral curve \eqref{kurva} given by the eigenvectors of the matrix $W(z)$. Let $D_0$ be the divisor of poles of a section of this line bundle and denote
\beq\label{tochka}
\bu_0=D_0-D_\infty-\Delta\in J(C)
\eeq
the point of the Jacobian corresponding to the line bundle. Here $D_\infty$ is the divisor of poles of the function $z:C\to {\bf P}^1$, $\Delta$ is the Riemann divisor (see \cite{Fay} for the definition). Here and below we identify divisors of degree 0 with their images in the Jacobian by means of the Abel--Jacobi map $A:C\to J(C)$,
$$
A(P-Q)=\left( \int_Q^P \w_1, \dots ,\int_Q^P \w_g\right).
$$

Our interest is in the $N$-differentials
\beq\label{N3}
\sum_{k_1, \dots, k_N=1}^g \frac{\pal^N\log\theta(\bu_0)}{\pal u_{k_1}\dots \pal u_{k_N}} \,\w_{k_1}(Q_1)\dots \w_{k_N}(Q_N), \quad Q_1, \dots, Q_N\in C
\eeq
for any $N\geq 3$. For $N=2$ instead of \eqref{N3} we will be looking at the following bi-differential
\beq\label{N2}
\frac{\theta\left( P-Q-\bu_0\right)\theta\left( P-Q+\bu_0\right)}{\theta^2(\bu_0) E(P,Q)^2}, \quad P, \, Q\in C
\eeq
that differs from the second logarithmic derivative of the form \eqref{N3} by the fundamental \emph{normalized} bi-differential (see Corollary 2.12 in the J.Fay's book \cite{Fay}). Here $E(P,Q)$ is the prime-form [{\it ibid}.]
$$
E(P,Q)=\frac{\theta[\nu](P-Q)}{\sqrt{\sum\omega_i(P)\pal_{u_i} \theta[\nu](0)}\sqrt{\sum\omega_j(Q) \pal_{u_j}\theta[\nu](0)}}
$$
where $\nu$ is a non-degenerate odd half-period. The goal is to express these multi-differentials on a spectral curve in terms of the associated matrix polynomial $W(z)$. The expression will involve the following matrix-valued function on the spectral curve
$C$
\beq\label{phi}
\Phi(P)\equiv \Phi(z,w)=\frac{R\left(z,W(z)\right)-R(z,w)}{W(z)-w}, \quad P=(z,w)\in C.
\eeq
In the right hand side it is understood that, first one has to cancel the  factor $W-w$ common for the numerator and denominator in the ratio $\frac{R(z,W)-R(z,w)}{W-w}$ and then to replace $W$ with $W(z)$. 

\medskip

{\bf Main Theorem}. {\it For the spectral curve \eqref{kurva} and for arbitrary points $Q_1=(z_1, w_1)$, \dots, $Q_N=(z_N, w_N)\in C$ the following expressions hold true
\beq\label{F2}
\frac{\theta\left( Q_1-Q_2-\bu_0\right)\theta\left( Q_1-Q_2+\bu_0\right)}{\theta^2(\bu_0) E(Q_1,Q_2)^2}=\tr \,\frac{\Phi(Q_1)\Phi(Q_2)}{(z_1-z_2)^2} \frac{dz_1}{R_w(z_1, w_1)} \frac{dz_2}{R_w(z_2, w_2)}
\eeq
\eqa\label{F3}
&&
\sum_{k_1, \dots, k_N=1}^g \frac{\pal^N\log\theta(\bu_0))}{\pal u_{k_1}\dots \pal u_{k_N}} \,\w_{k_1}(Q_1)\dots \w_{k_N}(Q_N)=
\\
&&
=\frac{(-1)^{N+1}}{N}\sum_{s\in S_N} \frac{\tr\, \left[\Phi\left( Q_{s_1}\right)\dots \Phi\left( Q_{s_N}\right)\right]}{{\left(z_{s_1}-z_{s_2}\right)\dots \left( z_{s_{N-1}} -z_{s_N}\right) \left( z_{s_N}-z_{s_1}\right)}} \frac{dz_1}{R_w(z_1,w_1)}\dots \frac{dz_N}{R_w(z_N,w_N)}
\nn
\eeqa
for any $N\geq 3$.
}

In these equations $R_w(z,w)$ is the partial derivative of the characteristic polynomial $R(z,w)$ with respect to $w$, the summation in \eqref{F3} is taken over all permutations of $\{ 1, 2, \dots, N\}$.

\medskip

For given $m$, $n$ denote $\CW_{m,n}$ the space of $n\times n$ matrix polynomials of the form \eqref{mapol}. Functions on $\CW_{m,n}$ will be denoted by $f\left([W]\right)$. By $\ZZ\left[ \CW_{m,n}\right]$
we denote the ring
$$
\ZZ\left[ \CW_{m,n}\right]=\ZZ\left[b^0_1, \dots, b^0_n, B^1_{ij},\dots, B^m_{ij},\prod_{i<j}(b^0_i-b^0_j)^{-1} \right].
$$
Let ${\mathcal C}_{m,n}$ be the space of algebraic curves of the form \eqref{kurva}, \eqref{nwcurva}. There is a natural fibration
\beq\label{fibra}
\begin{array}{c} \CW_{m,n}\\ \downarrow \\ {\mathcal C}_{m,n}\end{array}
\eeq
assigning to a matrix polynomial $W(z)$ its spectral curve $C$. The fiber over the point $C$ is isomorphic to the affine part of the generalized Jacobian (see below) of the 
singularized curve $C_{\rm sing}$ obtained from $C$ by identifying its infinite points, cf. \cite{DN74}, \cite{Mum}, \cite{NS}, \cite{SZ}, \cite{Gavril}. The quotient of the fiber over the diagonal conjugations
$$
W(z)\mapsto D^{-1}W(z) \, D, \quad D=\diag (d_1, \dots, d_n).
$$
is naturally isomorphic to $J(C)\setminus (\theta)$.

Define a function
\beq\label{Ft0}
F\left([W]; \bt\right)\in \ZZ\left[ \CW_{m,n}\right]\otimes \ZZ[[\bt]], \quad \bt= (t^a_k), \quad a=1, \dots, n, \quad k\geq 0
\eeq
by the infinite sum
\beq 
\label{Ft1}
F\left([W]; \bt\right)=\sum_{N=2}^\infty\frac1{N!}\sum_{a_1, \dots, a_N=1}^n\sum_{k_1,\dots, k_N=0}^\infty F_{k_1\dots k_N}^{a_1\dots a_N}[W] t^{a_1}_{k_1}\dots t^{a_N}_{k_N}
\eeq
where the coefficients $F_{k_1\dots k_N}^{a_1\dots a_N}[W]\in \QQ\left[ \CW_{m,n}\right]$ are defined by the following generated series
\beq\label{Ft2}
\sum_{k_1, k_2=0}^\infty \frac{F_{k_1 k_2}^{a_1 a_2}[W]}{z_1^{k_1}z_2^{k_2}}=\frac{\tr\,[\Pi_{a_1}(z_1)\Pi_{a_2}(z_2)]}{(z_1-z_2)^2}-\frac{\delta_{a_1, a_2}}{(z_1-z_2)^2}
\eeq
for any $a_1$, $a_2=1, \dots, n$ and
\beq\label{Ft3}
\sum_{k_1,\dots, k_N=0}^\infty \frac{F_{k_1\dots k_N}^{a_1\dots a_N}[W]}{z_1^{k_1+2}\dots z_N^{k_N+2}}=-\frac1{N}\sum_{s\in S_N} \frac{\tr\, \left[ \Pi_{s_1}(z_{s_1})\dots \Pi_{s_N}(z_{s_N})\right]}{{\left(z_{s_1}-z_{s_2}\right)\dots \left( z_{s_{N-1}} -z_{s_N}\right) \left( z_{s_N}-z_{s_1}\right)}}
\eeq
for any $N\geq 3$ and any $a_1$, \dots, $a_N=1,\dots, n$. Here the matrix-valued series 
$$
\Pi_a(z)=\Pi_a([W];z)\in Mat_n\left(\ZZ[\CW_{m,n}]\right)\otimes \ZZ[[z^{-1}]]
$$ 
are defined by the expansions
\beq\label{pia}
\frac{\Phi(P)}{R_w(z,w)}=\Pi_a(z), \quad P=(z,w)\to P_a,\quad a=1, \dots, n.
\eeq
Observe that the function $F\left([W]; \bt\right)$ is invariant with respect to diagonal conjugations.

{\bf Main Lemma}. {\it Let $W(z)\in \CW_{m,n}$ be any matrix-valued polynomial such that its spectral curve \eqref{kurva} has $n$ distinct points at infinity and it is nonsingular. Denote $\theta(\bu)$ the theta-function of the curve with respect to some basis of cycles and $\bu_0$ the point \eqref{tochka} of the Jacobian specified by the line bundle of eigenvectors of $W(z)$. Then the following equality of formal series in $\bt$ takes place
\beq\label{idprin}
e^{ F\left([W]; \bt\right)}=e^{\alpha +\sum \beta_{a,i} t^a_i +\frac12\sum \gamma_{a,i; b,j} t^a_i t^b_j}\,\theta\left(\sum t^a_k\bv^{(a,k)}-\bu_0\right)
\eeq
for some coefficients $\alpha$, $\beta_{a,i}$, $\gamma_{a,i; b,j}\in\mathbb C$. Here the vectors $\bv^{(a,k)}=\left(V^{(a,k)}_1, \dots, V^{(a,k)}_g\right)$ come from the coefficients of the expansion
\beq\label{nwholo}
\omega_i(P) =\sum_{k=0}^\infty \frac{V^{(a,k)}_i}{z^{k+2}} dz, \quad z=z(P), \quad P\to P_a
\eeq
of the holomorphic differentials $\omega_i(P)$ near infinity,
$i=1, \dots, g$, $a=1, \dots, n$.
}

\medskip


\begin{cor} \label{cor01} For any matrix-valued polynomial $W(z)\in \CW_{m,n}$ satisfying the conditions of Main Lemma the series \eqref{Ft1} has a non-zero radius of convergence when restricted onto a finite number of the indeterminates $\bt$.
\end{cor}

{\bf Definition}. {\it We say that $\det\left( w\cdot {\bf 1}-W(z)\right)=0$ is an algebraic spectral curve over $\QQ$ if $W(z)\in Mat_n(\QQ)\otimes \QQ[z]$.
}

\begin{cor} \label{cor02} For an algebraic spectral curve over $\QQ$ all coefficients of expansions of the multi-differentials \eqref{N3}, \eqref{N2} near infinity are rational numbers.
\end{cor}

\medskip

One can also derive from the Main Theorem  some identities between sums of products of Riemann theta-function and its logarithmic derivatives. 

\begin{cor}\label{cor03} Let $C$ be an arbitrary compact Riemann surface of genus $g>0$. For any $N\geq 3$ the following identity holds true for its theta-function along with the normalised holomorphic differentials on $C$
\eqa\label{id3}
&&
\theta^N(\bu)\sum_{i_1, \dots, i_N=1}^g\frac{\pal^N\log\theta(\bu)}{\pal u_{i_1}\dots \pal u_{i_N}} \omega_{i_1}(Q_1)\dots \omega_{i_N}(Q_N)=
\\
&&
=-\frac1{N} \sum_{s\in S_N}\frac{\theta(Q_{s_1}-Q_{s_2}+\bu) \theta(Q_{s_2}-Q_{s_3}+\bu)\dots \theta(Q_{s_{N-1}}-Q_{s_N}+\bu) \theta(Q_{s_N}-Q_{s_1}+\bu)}{E(Q_{s_1}, Q_{s_2})E(Q_{s_2}, Q_{s_3})\dots E(Q_{s_{N-1}}, Q_{s_N}) E(Q_{s_N}, Q_{s_1})}
\nn
\eeqa
for arbitrary points $Q_1$, $Q_2$, \dots, $Q_N\in C$ and arbitrary $\bu\in J(C)\setminus (\theta)$.
\end{cor}

\medskip

See below eqs. \eqref{ex3}, \eqref{ex4} for the explicit spelling of the above identity for the cases $N=3$ and $N=4$. For $N=4$ and $\bu=0$ the identity \eqref{id3} appeared in \cite{Fay} (see Proposition 2.14 there). We did not find  in the literature other identities of the form \eqref{id3}.

The constructions of the present paper can be extended to algebraic spectral curves with an arbitrary ramification profile at infinity. This will be done\footnote{For hyperelliptic curves with a branch point at infinity this has already been done in \cite{Dub18}.} a subsequent publication.

Before we proceed to the precise constructions and to the proofs let us say few words about the main ideas behind them. First, it is the connection between spectral curves and their theta-functions with particular classes of solutions to integrable systems, see e.g. \cite{Krich77, Dub81}. Second, it is the remarkable idea that goes back to M.Sato \emph{et al.} that suggests to consider tau-functions of integrable systems as partition functions of quantum field theories, see e.g. \cite{Date, Segal}. The time variables of the integrable hierarchy play the role of coupling constants. So, the logarithmic derivatives of tau-functions can be considered as the connected correlators of the underlined quantum field theory. For the algebro-geometric solutions the tau-function essentially coincides with the theta-function of the spectral curve, up to multiplication by exponential of a quadratic form. There were quite a few interesting results in the theory of theta-functions inspired by this connection, see e.g. \cite{Kawa, Naka}. The novelty of approach of the present work is that we are looking more on the correlators than on the tau-function. The main tools in proving the statements formulated above is in using the algorithm of \cite{BDY1, BDY2, BDY3} developed for efficient computation of correlators in cohomological field theories. This algorithm applied to tau-functions of algebraic spectral curves readily produces the explicit expressions for the correlators given above.

\setcounter{equation}{0}
\setcounter{theorem}{0}
\section{Main constructions and proofs}\par

\setcounter{equation}{0}
\setcounter{theorem}{0}
\subsection{Matrix polynomials $\leftrightarrow$ spectral curves + divisors}\label{sec21}\par

For a given $n\geq 2$, $m\geq 1$ for $n>2$ or $m\geq 2$ for $n=2$ consider the space $\CW$ of matrix polynomials of the form
\beq\label{spacenw}
\CW=\{ W(z) =z^m B +\mbox{lower degree terms}\}
\eeq
where $B=\diag (b_1, \dots, b_n)$ is an arbitrary diagonal matrix. For any $W(z)\in \CW$ the corresponding spectral curve $C$ is
of the form
\beq\label{nwcurva}
\det\left( w\cdot{\bf 1} -W(z)\right)=w^n+a_1(z) w^{n-1}+\dots+a_n(z)=0, \quad \deg a_i(z)=m\, i, \quad i=1, \dots, n.
\eeq
Assume the entries $b_1$, \dots, $b_n$ of the matrix $B$ 
to be pairwise distinct. Then the Riemann surface \eqref{nwcurva} has $n$ distinct points $P_1$, \dots, $P_n$ at infinity,
$$
P_a=\{ z\to\infty, ~ w\to\infty, ~ \frac{w}{z^m}\to b_a\}, \quad a=1, \dots, n.
$$
For the algebraic curve \eqref{nwcurva} this condition translates as follows. Let $a_i(z) =\alpha_i z^{m\, i}+\dots$, $i=1, \dots, n$. Then the roots of the equation
\beq\label{nwcurva10}
b^n+\alpha_1 b^{n-1}+\dots+\alpha_n=0
\eeq
must be pairwise distinct.

Assuming smoothness of the finite part of the curve we compute its genus
\beq\label{nwg}
g=\frac{(n-1)(m\, n-2)}2.
\eeq

We have a natural line bundle ${\mathcal L}$ over $C$ of the eigenvectors
\beq\label{line}
W(z) {\bpsi}(P)=w\,{\boldsymbol \psi}(P), \quad P=(z,w)\in C, \quad \bpsi(P)=(\psi^1(P), \dots, \psi^n(P))^T
\eeq
(the symbol $(\,.\,)^T$ stands for the transposition) of the matrix $W(z)$. We will associate with this line bundle a point in the \emph{generalized Jacobian} $J(C; P_1, \dots, P_n)$ that can be considered as an analogue of the Jacobi variety for the singular curve obtained by gluing together all infinite points $P_1$, \dots, $P_n$, see, e.g., \cite{Gavril}. It can be represented by classes of \emph{relative linear equivalence} of divisors of degree zero on the curve. By definition two divisors $D_1$ and $D_2$ of the same degree belong to the same relative linear equivalence class if there exists a rational function $f$ on the curve $C$ with $(f)=D_1-D_2$ satisfying $f(P_1)=f(P_2)=\dots =f(P_n)$. There is a natural fibration
\beq\label{fibra}
J(C; P_1, \dots, P_n)\to J(P)
\eeq
associating with any divisor its class of linear equivalence.

With the line bundle ${\mathcal L}$ we associate a divisor $D_0$  on the spectral curve defined by
\beq\label{D0}
D_0=\{ P\in C\, |\,\psi^1(P)+\dots+\psi^n(P)=0\}
\eeq
for a nonzero eigenvector. It can be considered as the divisor of poles of the section normalized by the condition 
\beq\label{nwnorm}
\psi^1+\dots+\psi^n=1.
\eeq
The components of the eigenvector can be represented as
\beq\label{psidelta}
\psi^j(P)=\frac{\Delta_{ij}(z,w)}{\sum_{s=1}^n \Delta_{is}(z,w)}, \quad P=(z,w)\in C, \quad j=1, \dots, n
\eeq
for any $i$. Here $\Delta_{ij}(z,w)$ is the $(i,j)$-cofactor of the matrix $w\cdot {\bf 1}-W(z)$. So at infinity the normalised \eqref{nwnorm} eigenvector behaves as
\beq\label{nwinf}
\psi^i(P)=\delta_{i\,a} +{\mathcal O}\left( \frac1{z(P)}\right), \quad P\to P_a.
\eeq


\begin{lemma} \label{lem211} The divisor $D_0$ of poles of the eigenvector normalised by eq. \eqref{nwnorm} is a nonspecial divisor on $C\setminus (P_1\cup\dots\cup P_n)$ of degree $g+n-1$. Conversely, for any nonsingular curve $C$ of the form 
\beq\label{nwcurva1}
w^n+a_1(z) w^{n-1}+\dots+a_n(z)=0, \quad \deg a_i(z)=m\, i, \quad i=1, \dots, n
\eeq
with $n$ distinct ordered points $P_1$, \dots, $P_n$ at infinity and an arbitrary nonspecial divisor $D_0\subset C\setminus (P_1\cup\dots\cup P_n)$ of degree $g+n-1$ there exists a unique matrix polynomial $W(z)$ of the form \eqref{spacenw} with the spectral curve coinciding with $C$ and the divisor of poles of the normalized eigenvector coinciding with $D_0$.
\end{lemma}

\pf To prove the first part of Lemma we will use an algorithm developed in \cite{dien} for computing the poles of eigenvectors of matrix polynomials adjusting it to the present situation. Denote $e^*=(1,1,\dots, 1)\in{\mathbb C^n}^*$ and define a polynomial
\beq\label{nwdz}
D(z)=e^*\wedge e^* W(z) \wedge \dots \wedge e^* W^{n-1}(z)\in\Lambda^n{\mathbb C^n}^*\otimes \mathbb C[z].
\eeq
We also define polynomials $q_{ij}(z)$, $i, \, j=1, \dots, n$ as coefficients of the expansion of sums of cofactors $\Delta_{ij}(z,w)$ of the matrix $w\cdot {\bf 1}-W(z)$
\beq\label{nwq}
\sum_{s=1}^n \Delta_{is}(z,w)=w^{n-1}+q_{i\,2}(z) w^{n-2}+\dots+q_{i, n-1}(z) w+ q_{i, n}(z), \quad i=1, \dots, n.
\eeq
($q_{i1}(z)=1$ for any $i=1, \dots, n$).

\begin{prop}\label{prop212} 1) For any matrix polynomial of the form
$$
W(z)=z^m B+{\rm lower~degree~terms}, \quad B=\diag(b_1, \dots, b_n), \quad b_i\neq b_j\quad{\rm for}\quad i\neq j
$$
the polynomial \eqref{nwdz} has degree $\frac{mn(n-1)}2=g+n-1$.

2) Assume that the spectral curve of $W(z)$ is nonsingular and the roots $z_1$, \dots, $z_{g+n-1}$ of the polynomial $D(z)$ are pairwise distinct. Then the rank of the rectangular matrix
$$
C(z)=\left( q_{ij}(z) \right)_{1\leq i\leq n, ~ 1\leq j\leq n-1}
$$
evaluated at $z=z_k$, $k=1, \dots, g+n-1$, is equal to $n-1$.

3) For a given $k\in\{1, 2, \dots, g+n-1\}$ let $C_k$ be a non-zero $(n-1)$-minor of the matrix $C(z_k)$,
$$
C_k =\left( q_{i_s ,j}(z_k)\right)_{1\leq s,\, j\leq n-1}.
$$
Denote $\hat C_k$ the matrix obtained from $C_k$ by changing the last column
$$
q_{i_s,n-1}(z_k) \mapsto q_{i_s, n}(z_k), \quad s=1, \dots, n-1
$$
and put
\beq\label{nww} 
w_k=-\frac{\det \hat C_k}{\det C_k}, \quad k=1, \dots, g+n-1.
\eeq
Then the poles of the eigenvector of the matrix $W(z)$ normalized by \eqref{nwnorm} are at the points
\beq\label{nwpole}
Q_1=(z_1, w_1), \dots, Q_{g+n-1}=(z_{g+n-1}, w_{g+n-1}) \in C\setminus \{ P_1\cup\dots P_n\}.
\eeq
\end{prop}

Let us prove that the divisor $D_0$ is nonspecial, i.e., that $l(D_0)=n$, where $l({D})={\rm dim}\, H^0\left(C, {\mathcal O}({D})\right)$. Indeed, if $l(D_0)>n$ then there exists a non-constant rational function $f$ on the curve $C$ with poles\footnote{Here and below we will say that a rational function $f$ on the curve $C$ has poles at the points of a divisor $D$ if $(f)+D\geq 0$.}  at $D_0$ satisfying $f(P_1)=\dots=f(P_n)=1$. Consider the vector-function
$$
\tilde\bpsi(P)=\frac1{f(P)} \bpsi(P).
$$
Clearly it is again an eigenvector of the matrix $W(z)$ satisfying the normalization \eqref{nwnorm}. Due to uniqueness it must coincide with $\bpsi$, so $f$ must be identically equal to 1. Such a contradiction completes the proof of the first part of Lemma. 

Let us now explain the reconstruction procedure of the polynomial matrix $W(z)$ starting from a pair $(C,D_0)$ consisting of a curve $C$ of the form \eqref{nwcurva} smooth for $|z|<\infty$ and a nonspecial positive divisor $D_0$ on $C\setminus (P_1\cup\dots\cup P_n)$ of degree $g+n-1$. As by assumption $l(D_0)=n$, there exist $n$ rational functions $\psi^1(P)$, \dots, $\psi^n(P)$ on $C$ with poles at $D_0$ satisfying
\beq\label{normpsi}
\psi^i(P_j)=\delta_{ij}, \quad i, \, j=1, \dots, n.
\eeq
Let $R$ be a sufficiently large number such that no ramification neither points of the divisor $D_0$ occur for $|z|>R$. For any such $z$
denote $(z,1)$, \dots, $(z,n)$ the preimages of $z$ on the spectral curve with respect to the natural projection $z:C\mapsto {\bf P}^1$ ordered in an arbitrary way. Define a $n\times n$ matrix $\Psi(z)$ whose $i$-th row is $\left( \psi^i\left((z,1)\right), \dots, \psi^i\left((z,n)\right)\right)$. The matrix $\Psi(z)$ is invertible for any $R<|z|\leq\infty$. Denote
$$
w_a(z)=w\left((z,a)\right)=b_a z^m+\dots, a=1, \dots, n
$$
the branches of the algebraic function $\mu(P)$, $P\in C$ and put
$$
\hat w(z)=\diag (w_1(z),\dots, w_n(z)),\quad  |z|>R.
$$
The matrix-valued function $\Psi(z)\hat w(z)\Psi^{-1}(z)$ is analytic for $|z|>R$ having an $m$-th order pole at infinity. Observe that it does not depend on the ordering of the preimages of $z$, so it can be extended to a rational function on the complex plane. Consider its Laurent expansion at infinity
$$
\Psi(z)\hat w(z)\Psi^{-1}(z)=B z^m +B_1z^{m-1}+\dots+B_m+{\mathcal O}\left(\frac1{z}\right)
$$
where $B=\diag (b_1,\dots, b_n)$ and $B_1$,\dots, $B_m$ are some $n\times n$ matrices. Put
$$
W(z)=B z^m +B_1z^{m-1}+\dots+B_m
$$
and consider the vector-valued function $\tilde\psi(P)$ on the curve defined by
$$
\tilde\psi(P)=W\left(z(P)\right) \psi(P)-w(P) \psi(P).
$$
It has poles only at the points of the divisor  $D_0$. From the definition of the matrix $W(z)$ it readily follows that $\tilde\psi(P)$ vanishes at $P_1$,\dots, $P_n$. Hence it equals zero due to the nonspeciality of the divisor $D_0$. This implies that $C$ coincides with the spectral curve of the matrix polynomial $W(z)$. To complete the reconstruction procedure it remains to observe that the function
$$
\psi^1(P)+\dots+\psi^n(P)-1
$$
having poles at $D_0$
vanishes at $P_1$, \dots, $P_n$. Hence it is identically equal to 0, that is, the eigenvector $\psi(P)$ of the matrix $W(z)$ satisfies the normalization \eqref{nwnorm}. \epf

\begin{remark} Changing the divisor $D_0\to D_0'\sim D_0$ in the class of linear equivalence yields a conjugation of the matrix $W(z)$ by a diagonal matrix
$$
W(z) \to F\, W(z) F^{-1}, \quad F=\diag\left(f(P_1), \dots, f(P_n)\right)
$$
where $f$ is a rational function on $C$ with the divisor $D_0-D_0'$.
\end{remark}

One can repeat the above constructions dealing with the \emph{dual line bundle} ${\mathcal L}^\dagger$ over $C$ coming from the left eigenvectors of the matrix $W(z)$
\beq\label{duline}
{\bpsi}^\dagger(P)\,W(z) =w\,{ \bpsi}^\dagger(P), \quad P=(z,w)\in C, \quad \bpsi^\dagger(P)=(\psi_1^\dagger(P), \dots, \psi^\dagger_n(P)).
\eeq
Let us use the same normalization
\beq\label{dunwnorm}
\psi_1^\dagger+\dots +\psi_n^\dagger=1
\eeq
so
\beq\label{dupsidelta}
\psi_i^\dagger(P)=\frac{\Delta_{ij}(z,w)}{\sum_{s=1}^n\Delta_{sj}(z,w)},\quad P=(z,w)\in C, \quad i=1, \dots n
\eeq
for any choice of $j$ (cf. eq. \eqref{psidelta}). Denote $D_0^\dagger$ the divisor of poles of the dual eigenvector \eqref{duline} normalized by the condition \eqref{dunwnorm}.
The following statement is an analogue of Lemma \ref{lem211} for the dual eigenvectors.

\begin{lemma}\label{lem2110du} The divisor $D_0^\dagger$ of poles of the eigenvector \eqref{duline} normalised by eq. \eqref{dunwnorm} is a nonspecial divisor on $C\setminus (P_1\cup\dots\cup P_n)$ of degree $g+n-1$. Conversely, for any nonsingular curve $C$ of the form 
\beq\label{nwcurva11}
w^n+a_1(z) w^{n-1}+\dots+a_n(z)=0, \quad \deg a_i(z)=m\, i, \quad i=1, \dots, n
\eeq
and an arbitrary nonspecial divisor $D_0^\dagger\subset C\setminus (P_1\cup\dots\cup P_n)$ of degree $g+n-1$ there exists a unique matrix polynomial $W(z)$ of the form \eqref{spacenw} with the spectral curve coinciding with $C$ and the divisor of the normalized eigenvector \eqref{duline} coinciding with $D_0$.
\end{lemma}

The proof is similar to that of Lemma \ref{lem211}, so it will be omitted.

\begin{remark} Of course the divisors $D_0$ and $D_0^\dagger$ do depend on each other. The nature of this dependence will be clarified below.
\end{remark}

We will now look at the \emph{spectral projectors} of the matrix polynomial $W(z)$.
Consider the matrix-valued function 
\beq\label{nwpi0}
\Pi(P)=\frac{\Phi(P)}{R_w(z,w)},\quad P=(z,w)\in C
\eeq
where $R(z,w)=\det\left(w\cdot{\bf 1}-W(z)\right)$ is the characteristic polynomial of the matrix $W(z)$ and $\Phi(P)$ is defined by
\eqref{phi}. So
\eqa\label{piex}
&&
\Pi(z,w)=\frac1{R_w(z,w)}\frac{R(z,W)-R(z,w)}{W-w} =\frac{\sum_{i=0}^{n-1}b_i(z) w^{n-i-1}}{R_w(z,w)}
\\
&&
b_i(z)=\sum_{j=0}^i a_j(z) W^{i-j}(z).
\nn
\eeqa
Let $z\in\mathbb C$ be not a ramification point. Denote $w_1(z)$, \dots, $w_n(z)$ the points above $z$ on the spectral curve $C$ and put
\beq\label{pii}
\Pi_i(z)=\Pi(z,w_i(z)), \quad i=1, \dots, n.
\eeq
We will prove that these matrices are the spectral projectors
$$
\Pi_i(z): \mathbb C^n \to {\rm Ker}\,\left( W(z)-w_i(z)\cdot {\bf 1}\right), \quad i=1, \dots, n
$$
of $W(z)$.

\begin{lemma} The matrices $\Pi_1(z)$, \dots, $\Pi_n(z)$ are basic idempotents of the matrix $W(z)$, i.e.
\eqa\label{id}
&&
\Pi_i^2=\Pi_i, \quad \Pi_i\cdot \Pi_j=0\quad{\rm for}\quad i\neq j, \quad i, \, j=1, \dots, n
\\
&&\label{id1}
\sum_{i=1}^n \Pi_i(z)={\bf 1}, \quad \sum_{i=1}^n w_i(z) \Pi_i(z)=W(z).
\eeqa
\end{lemma}

\pf Let us first prove that
\beq\label{syst}
\sum_{i=1}^n w_i^r(z) \Pi_i(z)=W^r(z) \quad \mbox{for any}\quad 0\leq r.
\eeq
We have
$$
\sum_{i=1}^n w_i^r(z)\frac{\sum_{k=0}^{n-1}b_k(z) w_i^{n-k-1}(z)}{R_w(z,w_i(z))}=\sum_{i=1}^n \res_{w=w_i(z)} w^r\frac{\sum_{k=0}^{n-1}b_k(z) w^{n-k-1}}{R(z,w)} =-\res_{w=\infty} w^r\frac{\sum_{k=0}^{n-1}b_k(z) w^{n-k-1}}{R(z,w)}.
$$
Using the explicit expression \eqref{piex} along with the obvious identity
$$
\left(1-\frac{W}{w}\right)\left(1+\frac{b_1}{w}+\dots+\frac{b_{n-1}}{w^{n-1}}\right)=1+\frac{a_1}{w}+\dots +\frac{a_n}{w^n}
$$
we arrive at
\beq\label{expan}
\frac{\sum_{k=0}^{n-1}b_k(z) w^{n-k-1}}{R(z,w)}=\frac1{w-W}, \quad |w|\to\infty.
\eeq
Thus
$$
\res_{w=\infty} w^r\frac{\sum_{k=0}^{n-1}b_k(z) w^{n-k-1}}{R(z,w)}=-W^r.
$$
This proves \eqref{syst} and, hence \eqref{id1}.

To prove \eqref{id} we solve the system \eqref{syst} for $r=0$, $1$, \dots, $n-1$ with respect to $\Pi_1(z)$, \dots $\Pi_n(z)$ to obtain
$$
\Pi_i(z)=\frac{\prod_{j\neq i} (W(z)-w_j(z))}{\prod_{j\neq i} (w_i(z)-w_j(z))}, \quad i=1, \dots, n.
$$
Rewriting these matrices in the basis of eigenvectors of $W(z)$ we readily get \eqref{id}. \epf

It will be convenient to also consider a matrix-valued differential with the matrix entries
\beq\label{matdif}
\Omega^i_j(P)=\Pi^i_j(P)dz, \quad i,\, j=1, \dots, n
\eeq
where the matrix $\Pi(P)=\left( \Pi^i_j(P)\right)_{1\leq i,\, j\leq n}$ is given by \eqref{nwpi0}, \eqref{piex}.

\begin{prop} \label{prop217} For every $i\neq j$ the differential $\Omega^i_{j}(P)$ is holomorphic on $C\setminus \left(P_i \cup P_j\right)$ with simple poles at $P=P_i$ and $P=P_j$.
The differential $\Omega^i_{i}(P)$ is holomorphic on $C\setminus P_i$ having a double pole at $P=P_i$ such that
\beq\label{matdif1}
\Omega^i_i(P)={dz}+{\rm regular~terms}, \quad P\to P_i
\eeq
for every $i=1, \dots, n$.
\end{prop}

\pf As it follows from the explicit expression \eqref{phi} the entries of the matrix $\Phi(P)$ are holomorphic on the affine part $P=(z,w)\in C$, $|z|<\infty$. The differential $dz/R_w(z,w)$ is holomorphic on $C$. Hence the differentials $\Pi_{ij}(P)$ are holomorphic on the affine part of $C$. Let us look at their behaviour at infinity. To this end we will use the standard realization of the spectral projectors of a $n\times n$ complex matrix with pairwise distinct eigenvalues that in our case can be formulated in the following way. Let $z\in\mathbb C$ be not a ramification point with respect to the projection $z: C\to \mathbb C$. Order the points $(z,w_1(z))$, \dots, $(z,w_n(z))$ in the preimage. Like in the proof of Lemma \ref{lem211} produce a matrix $\Psi(z)$ whose $k$-th column is given by the eigenvector $\left(\psi^1(z,w_k(z)), \dots, \psi^n(z,w_k(z))\right)^T$ of $W(z)$ normalized by the condition \eqref{nwnorm}. Then
\beq\label{basic}
\Pi(P)|_{P=(z,w_k(z))}=\Psi(z) E_k \Psi^{-1}(z)
\eeq
where the $n\times n$ matrix $E_k$ has only one nonzero entry
$$
\left(E_k\right)^i_j=\delta^i_k\delta^k_j.
$$
Consider now sufficiently large $R$; choose the order of eigenvalues over the disk $|z|>R$ in such a way that $w_k(z)\sim b_k z^m$ for $|z|\to \infty$. Then
\beq\label{beha}
\Psi(z)={\bf 1}+{\mathcal O}\left(\frac1{z}\right), \quad \Psi^{-1}(z)={\bf 1}+{\mathcal O}\left(\frac1{z}\right), \quad z\to\infty
\eeq
due to eq. \eqref{nwinf}. The behaviour of the differentials $\Omega^i_j(P)$ at infinity easily follows from \eqref{basic}, \eqref{beha}. \epf

Consider now the differentials
\beq\label{omlow}
\Omega_j(P)=\sum_{i=1}^n \Omega^i_j(P)
\eeq
\beq\label{omup}
\Omega^i(P)=\sum_{j=1}^n\Omega^i_j(P).
\eeq

\begin{lemma} \label{lem218} The differential $\Omega_j(P)$ has zeros at the points of the divisor $D$ and at some divisor $D^\dagger_j$ of degree $g$. It has poles at the points of the divisor $P_j+\sum_{s=1}^n P_s$. In a similar way the differential $\Omega^i(P)$ has zeros at the points of the divisor $D^\dagger$ and at some divisor $D_i$ of degree $g$.  It has poles at the points of the divisor $P_i+\sum_{s=1}^n P_s$. 
\end{lemma}

\pf Using the representation \eqref{basic} along with the normalization \eqref{nwnorm} we immediately conclude that the sum \eqref{omlow} vanishes at the points of the divisor $D$. The configuration of poles of this differential can be easily recovered from Proposition \ref{prop217}. The degree counting
$$
\deg D+\deg D_j^\dagger-\deg\left( P_j+\sum_{s=1}^n P_s\right)=2g-2
$$
yields $\deg D_j^\dagger=g$. To derive similar statements about the differentials $\Omega^i(P)$ we use an alternative representation of the projector matrix
\beq\label{dubasic}
\Pi(P)|_{P=(z,w_k(z))}={\Psi^\dagger(z)}^{-1} E_k \Psi^\dagger(z)
\eeq
where the $k$-th row of the matrix $\Psi^\dagger(z)$ is given by the left eigenvector $\left(\psi^\dagger_1(z,w_k(z)),\dots, \psi^\dagger_n(z,w_k(z))\right)$ normalized by the condition \eqref{dunwnorm}. \epf

\begin{cor}\label{cor219}  The divisor of zeros of the differential $\Omega^i_j(P)$ on $C\setminus \left( P_1\cup\dots P_n\right)$ coincides with $D_i+D_j^\dagger$, $i$, $j=1,\dots, n$.
\end{cor}

\pf The matrix $\Omega^i_j(P)$ has rank one. Its columns are eigenvectors of the matrix $W(z)$ with the same eigenvalue. Normalizing any column we obtain the same vector function $\bpsi(P)$
\beq\label{drob1}
\psi^i(P)=\frac{\Omega^i_j(P)}{\sum_{k=1}^n \Omega^k_j(P)}, \quad i=1,\dots , n
\eeq
for any $j$. According to Lemma the denominator vanishes at the points of the divisor $D_j^\dagger$. Hence also the numerator must vanish at the points of this divisor. In a similar way, the rows of $\Omega^i_j(P)$ are left eigenvectors of the same matrix $W(z)$. Normalizing them one obtains
\beq\label{drob2}
\psi^\dagger_j(P)=\frac{\Omega^i_j(P)}{\sum_{k=1}^n \Omega^i_k(P)}, \quad j=1,\dots , n
\eeq
for any $i$. Hence $\Omega^i_j(P)$ vanishes also at the points of the divisor $D_i$. Since degree of the divisor of poles of this differentials equals two and $\deg (D_i+D_j^\dagger)=2g$, there are no other zeros. \epf

\begin{cor} The zeros of the differential
\beq\label{difom1}
\Omega(P)=\sum_{i,\,j=1}^n\Omega^i_j(P)
\eeq
are at the points of the divisor $D+D^\dagger$. It has double poles at the infinite points $P_1$, \dots, $P_n$ and
\beq\label{difom2}
\Omega(P)=\left(1+{\mathcal O}\left(\frac1{z}\right)\right) dz, \quad P\to P_k, \quad k=1,\dots,n.
\eeq
\end{cor}

The corollary suggests the following way of determining the dual divisor $D^\dagger$ starting from $D$. For a nonspecial divisor $D$ of degree $g+n-1$ there exists a unique differential $\Omega(P)$ vanishing at the points of $D$ and having double poles of the form \eqref{difom2} at infinity. The remaining zeros of the differential give the points of the divisor $D^\dagger$.

\begin{cor} The differentials $\Omega^i_j(P)$ admit the following representation
\beq\label{omij}
\Omega^i_j(P)=\psi^i(P)\psi^\dagger_j(P) \Omega(P),\quad i,\, j=1,\dots, n.
\eeq
\end{cor}

\pf Due to the previous Corollary the product \eqref{omij} is holomorphic on $C\setminus\left( P_1\cup\dots\cup P_n\right)$. From Corollary \ref{cor219} it follows that the divisor of zeroes of this product coincides with $D_i+D^\dagger_j$. Finally, using \eqref{difom2} along with the asymptotics of $\psi^i(P)$ and $\psi^\dagger_j(P)$ at infinity we complete the proof. \epf

\setcounter{theorem}{0}
\subsection{Generalized Jacobian and theta-functions. Proof of eq. \eqref{F2}}\par

The generalized Jacobian $J(C; P_1, \dots, P_n)$ can be realized as a fiber bundle over $J(C)$ with $(n-1)$-dimensional fiber. For $n=2$ the construction was already given in \cite{Fay}; it is quite similar also for arbitrary $n$. 

Define $J(C; P_1, \dots, P_n)$ as the set of all pairs
$$
(\bu, \bla), \quad \bu=(u_1, \dots, u_g)\in\mathbb C^g, \quad \bla=\left( \lambda_1, \dots, \lambda_n\right)\in \left( \mathbb C^*\right)^n
$$
modulo the following equivalence relation $(\bu, \bla)\sim  (\bu', \bla')$ if
\beq\label{genjac}
\bu'=\bu+2\pi i M+BN,  \quad \lambda_k'=c\,\lambda_k \,e^{\langle N, A(P_k)\rangle}, \quad k=1, \dots, n,  \quad c\in \mathbb C^*, \quad M, \, N \in \mathbb Z^g.
\eeq
The fibration \eqref{fibra} is realized by the map $(\bu, \bla) \mapsto \bu$.

To define an analogue of the Abel map
\beq\label{abjac}
C\to J(C; P_1, \dots, P_n)
\eeq
fix a pair of distinct points $P_0$, $Q_0\in C\setminus \left( P_1\cup \dots \cup P_n\right)$ and put
\beq\label{abjac1}
P\mapsto \left( \bu(P), \bla(P)\right), \quad u_i(P)=A_i(P)=\int_{P_0}^P \omega_i, \quad \lambda_k(P) =e^{\alpha_k(P)}, \quad \alpha_k(P)= \int_{P_0}^P \Omega_{P_k Q_0}.
\eeq
Here and below $\Omega_{PQ}$ is the third kind differential on $C$ having simple poles at the points $P$ and $Q$ with residues $+1$ and $-1$ respectively and vanishing $a$-periods. The map is extended linearly/multiplicatively on the group of divisors of a given degree. The following statement is an analogue of the Abel--Jacobi theorem.

\begin{prop} Two divisors $D$, $D'$ of the same degree are relatively equivalent \emph{iff}
$$
\left(\bu(D), \bla(D)\right)\sim \left(\bu(D'), \bla(D')\right)
$$
modulo equivalence \eqref{genjac}.
\end{prop}

An analogue of the Riemann theorem about zeros of theta-function is given by

\begin{prop} \label{prop221} For a given $\bu\in\mathbb C^g$ such that $\theta(\bu)\neq 0$ and $\bla\in\left(\mathbb C^*\right)^n$ consider the function
\beq\label{fp}
F(P)=\sum_{s=1}^n \lambda_s \,\theta(P-P_s-\bu) \frac{E(P,P_0)(d\zeta_0)^{1/2}}{E(P,P_s)(d\zeta_s)^{1/2}}
\eeq
on the $4g$-gon $\tilde C$ obtained by cutting the curve $C$ along the chosen basis of cycles $a_1$, \dots, $a_g$, $b_1$, \dots, $b_g$. Here $\zeta_0(P)=z(P)-z(P_0)$, $\zeta_s(P)=1/z(P)$ are local parameters\footnote{For simplicity we assume that $P_0\in C$ is not a branch point.}  near $P_0$, $P_s$ respectively. This function has simple poles at $P_1$, \dots, $P_n$, a simple zero at $P_0$  and also zeros at some points $Q_1$, \dots, $Q_{g+n-1}$ of a divisor $D\subset C\setminus \left( P_1\cup \dots \cup P_n\right)$ satisfying
\eqa\label{riem1}
&&
\sum_{i=1}^{g+n-1} A_j\left( Q_i\right) -\sum_{s=1}^n A_j\left(P_s\right) =u_j-{\mathcal K}_j, \quad j=1, \dots, g
\\
&&
\label{riem2}
\exp \sum_{i=1}^{g+n-1} 
\alpha_s\left(Q_i\right)
=\lambda_s\, e^{-\kappa_s}, \quad s=1, \dots, n
\eeqa
where ${\mathcal K}_1$, \dots, ${\mathcal K}_g$ are Riemann constants
$$
{\mathcal K}_j=\frac{2\pi i+B_{jj}}2-\frac1{2\pi i}\sum_{k\neq j} \int_{a_k} A_j(P) \omega_k(P),
$$
$\kappa_1$, \dots, $\kappa_n$ are given by analogous formulae,
\beq\label{kappa}
\kappa_s=-\log E\left( P_s, P_0\right)-\frac1{2\pi i}\sum_{k=1}^g\int_{a_k}\alpha_s(P) \omega_k(P).
\eeq
\end{prop}

Recall \cite{Fay} that the combination
\beq\label{Delta}
\Delta:=-{\mathcal K}+(g-1)P_0
\eeq
does not depend on the choice of the base point $P_0$. This gives rise to definition of the Riemann divisor $\Delta$ that already appeared above. In particular eq. \eqref{riem1} can be rewritten in the form
\beq\label{riem11}
D -\sum_{s=1}^n P_s -\Delta= \bu.
\eeq
Observe that, for a given $(\bu, \bla )$ the divisor of zeros of the function \eqref{fp} determined by eqs. \eqref{riem1}, \eqref{riem2} does not depend on the choice of the basepoint $P_0$.

\pf The function $F(P)$ does not change its value when the point $P$ crosses the cut along the cycle $b_k$ while crossing the cut $a_k$ it is multiplied by
$$
\exp\left( -\frac12 B_{kk} -\int_{P_0}^P \omega_k +u_k\right).
$$
So the total logarithmic residue
$$
\frac1{2\pi \,i}\oint_{\pal\tilde C} d\log F(P)
$$
is equal to $g$. Hence $F(P)$ has $g+n-1$ zeros on $C\setminus \left( P_1\cup \dots \cup P_n\right)$, counted with multiplicity. To prove the first equality \eqref{riem1} one has to compute the contour integral
$$
\frac1{2\pi\, i} \oint_{\pal\tilde C} A_j(P) d\log F(P)=A_j(Q_1)+\dots+A_j(Q_{g+n-1})-A_j(P_1)-\dots-A_j(P_n).
$$
To prove eq. \eqref{riem2} we need to add more cuts: a cut from $P_0$ to $Q_0$ and also cuts from $Q_0$ to the infinite points $P_1$, \dots, $P_n$. Denote $\tilde C'$ the resulting polygon. Then
$$
\sum_{i=1}^{g+n-1}\alpha_s\left(P_i\right)
=\frac1{2\pi \, i} \oint_{\pal\tilde C'} \alpha_s(P)
d\log F(P), \quad s=1, \dots, n
$$
up to a $s$-independent shift.
The integral in the right hand side must be regularized at $P\to \infty$ or $P\to Q_0$. After such regularization we arrive at eqs. \eqref{riem2} up to equivalence \eqref{genjac}. \epf

We will now express the differentials \eqref{matdif} in terms of the coordinates $(\bu_0, \bla_0)$.

\begin{prop} \label{prop223} Let $Q_1+\dots+Q_{g+n-1}=D\subset C\setminus \left(P_1\cup\dots\cup P_n\right)$ be the divisor of poles of the eigenvector of the matrix $W(z)$ normalised by the condition \eqref{nwnorm}. Denote $(\bu_0, \bla_0)$ the corresponding point on the generalized Jacobian \eqref{genjac} 
\eqa\label{corp}
&&
\bu_0=D_0-\sum_{a=1}^n P_a -\Delta
\\
&&
\lambda_j^0=\exp\left\{ \sum_{s=1}^{g+n-1} \int_{P_0}^{Q_s} \Omega_{P_jQ_0}+\kappa_j\right\}, \quad j=1, \dots , n
\nn
\eeqa
(cf. eqs. \eqref{riem1}, \eqref{riem2}). Then the differentials $\Omega^i_j(P)$ of the form \eqref{matdif} are given by the following equation
\beq\label{glav}
\Omega^i_j(P)=\frac{\lambda_i^0}{\lambda_j^0}\frac{\theta(P-P_i-\bu_0)\theta(P-P_j+\bu_0)}{\theta^2(\bu_0)E(P_i,P)E(P,P_j)\sqrt{d\zeta_i}\sqrt{d\zeta_j}}
\eeq
\end{prop}

\pf Any differential $\Omega^i_j(P)$ having, for $i\neq j$ simple poles at $P=P_i$ and $P=P_j$ and, for $i=j$ a double pole of the form \eqref{matdif1} at $P=P_i$ can be written \cite{Fay} as follows
$$
\Omega^i_j(P)=\alpha_{ij} \frac{\theta(P-P_i-\bu_{ij})\theta(P-P_j+\bu_{ij})}{\theta^2(\bu_{ij})E(P_i,P)E(P,P_j)\sqrt{d\zeta_i}\sqrt{d\zeta_j}}
$$
for some $\bu_{ij}\in J(C)\setminus\left(\theta\right)$ and some nonzero constants $\alpha_{ij}$ satisfying $\alpha_{ii}=1$. According to Riemann theorem the zeros of the function $\theta(P-P_j+\bu_{ij})$ are at the points of a divisor ${\mathcal D}$ of degree $g$ satisfying
$$
{\mathcal D}-P_j-\Delta=-\bu_{ij}.
$$
According to the Corollary \ref{cor219} it must coincide with the divisor $D^\dagger_j$. From Lemma \ref{lem218} we derive the following linear equivalence
$$
D_0+D^\dagger_j-P_j-\sum_{a=1}^nP_a=K_C.
$$
Substituting $K_C=2\Delta$ we can rewrite it as follows
$$
D^\dagger_j-P_j-\Delta=-D_0+\sum_{a=1}^nP_a+\Delta=-\bu_0.
$$
Hence the condition ${\mathcal D}=D^\dagger_j$ implies $\bu_{ij}=\bu_0$ on $Jac(C)$.

It remains to fix the constants $\alpha_{ij}$. Since the rank of the matrix $\Omega^i_j(P)$ must be equal to one, we conclude that $\alpha_{ij}=\alpha_i \beta_j$ for some nonzero constants $\alpha_i$, $\beta_j$. As $\alpha_{ii}=1$ then $\beta_j=\alpha_j^{-1}$. The last condition to be used is that the sum $\sum_{i=1}^n\Omega^i_j(P)$ must vanish at the points $Q_1$, \dots, $Q_{g+n-1}$ of the divisor $D$. This implies that $\alpha_i=\lambda_i^0$, up to a common factor. \epf

\begin{cor} The eigenvector $\bpsi(P)$ of the matrix $W(z)$ normalized by the condition \eqref{nwnorm} is
\beq\label{nwpsi0}
\psi^i(P)=\frac{
\lambda_i^0 
\frac{
\theta(P-P_i-\bu_0)
}
{
 E(P,P_i) \left( d\zeta_i\right)^{1/2}}}
{\sum_{b=1}^n \lambda_b^0 \frac{\theta(P-P_b-\bu_0)}{E(P,P_b) \left( d\zeta_b\right)^{1/2}}},\quad i=1, \dots, n.
\eeq
The dual eigenvector $\bpsi^\dagger(P)$ is given by a similar formula
\beq\label{psidag}
\psi^\dagger_i(P)=\frac{
\frac1{\lambda_i^0 }
\frac{
\theta(P-P_i+\bu_0)
}
{
 E(P_i,P) \left( d\zeta_i\right)^{1/2}}}
{\sum_{b=1}^n \frac1{\lambda_b^0} \frac{\theta(P-P_b+\bu_0)}{E(P_b,P) \left( d\zeta_b\right)^{1/2}}},\quad i=1, \dots, n.
\eeq
\end{cor}

\pf Use \eqref{glav} along with \eqref{drob1}, \eqref{drob2}. \epf

\begin{remark} \label{rem225} Observe that the change
$$
D_0\mapsto D^\dagger_0
$$
corresponds to the involution
\beq\label{invo}
(\bu_0, \lambda^0_1, \dots,\lambda^0_n)\mapsto \left(-\bu_0, 1/\lambda_1^0,\dots, 1/\lambda_n^0\right)
\eeq
on the generalized Jacobian.
\end{remark}

We are now in a position to prove the first equation \eqref{F2} of the Main Theorem. 

\begin{prop} \label{prop226} Let $C$ be a compact Riemann surface of positive genus and $z:C\to {\bf P}^1$ a rational function with $n$ simple poles at the points $P_1$, \dots, $P_n$. Introduce the following matrix of Abelian differentials on $C$
\beq\label{matom}
\bom(P)=\left(\Omega^i_j(P)\right)_{1\leq i,\, j\leq n}, \quad \Omega^i_j(P)=\frac{\lambda_i}{\lambda_j}\frac{\theta(P-P_i-\bu)\theta(P-P_j+\bu)}{\theta^2(\bu) E(P_i,P) E(P,P_j)\sqrt{d\zeta_i}\sqrt{d\zeta_j}}
\eeq
where $\bu\in J(C)\setminus\left(\theta\right)$ is an arbitrary point and $\lambda_1$, \dots, $\lambda_n$ are arbitrary nonzero numbers. Then for an arbitrary pair of distinct points $P$, $Q\in C$ the following equation holds true
\beq\label{pref2}
\tr \frac{\bom(P)\bom(Q)}{\left(z(P)-z(Q)\right)^2}=\frac{\theta(P-Q-\bu)\theta(P-Q+\bu)}{\theta^2(\bu) E^2(P,Q)}.
\eeq
\end{prop}

\pf  The trace of the product of the matrices $\bom(P)$ and $\bom(Q)$ factorizes as follows
\beq\label{fac83}
\tr\,\bom(P)\bom(Q)
=\sum_{i=1}^n\frac{\theta(P-P_i-\bu)\theta(Q-P_i+\bu)}{\theta^2(\bu) E(P_i, P) E(Q, P_i) d\zeta_i} \cdot \sum_{j=1}^n\frac{\theta(Q-P_j-\bu)\theta(P-P_j+\bu)}{\theta^2(\bu) E(P_j, Q) E(P, P_j) d\zeta_j} .
\eeq
For a fixed pair of distinct points $P$, $Q\in C$ consider the differential
\beq\label{hpq}
H_{PQ}(Z)=\frac{\theta(P-Z-\bu)\theta(Q-Z+\bu)}{\theta^2(\bu)E(Z,P)E(Q,Z) \sqrt{dz(P)}\sqrt{dz(Q)}}, \quad Z\in C.
\eeq
It has simple poles at $Z=P$ and $Z=Q$ with residues
$$
\res_{Z=P} H_{PQ}(Z)=-\res_{Z=Q}H_{PQ}(Z)=\frac{\theta(P-Q-\bu)}{\theta(\bu)E(Q,P) \sqrt{dz(P)}\sqrt{dz(Q)}}.
$$
We now consider the product $z(Z) H_{PQ}(Z)$. Vanishing of the sum of residues of this differential yields
\eqa\label{sumid}
&&
\sum_{i=1}^n\frac{\theta(P-P_i-\bu)\theta(Q-P_i+\bu)}{\theta^2(\bu) E(P_i, P) E(Q, P_i)  \sqrt{dz(P)}\sqrt{dz(Q)}\,d\zeta_i}=\sum_{i=1}^n \res_{Z=P_i} z(Z)\, H_{PQ}(Z) =
\nn\\
&&
=\left[ z(P)-z(Q)\right] \frac{\theta(P-Q-\bu)}{\theta(\bu) E(Q,P) \sqrt{dz(P)}\sqrt{dz(Q)}}
\nn\\
&&
\\
&&
\sum_{j=1}^n\frac{\theta(Q-P_j-\bu)\theta(P-P_j+\bu)}{\theta^2(\bu) E(P_j, Q) E(P, P_j)\sqrt{dz(P)}\sqrt{dz(Q)}\, d\zeta_j}=\sum_{j=1}^n \res_{Z=P_j} z(Z)\, H_{QP}(Z)=
\nn\\
&&
=[z(P)-z(Q)]\frac{\theta(P-Q+\bu)}{\theta(\bu) E(Q,P) \sqrt{dz(P)}.\sqrt{dz(Q)}}
\nn
\eeqa
Therefore the bi-differential in the right hand side of eq. \eqref{fac83} becomes
\eqa
&&
\sum_{i=1}^n\frac{\theta(P-P_i-\bu)\theta(Q-P_i+\bu)}{\theta^2(\bu) E(P_i, P) E(Q, P_i) d\zeta_i} \cdot \sum_{j=1}^n\frac{\theta(Q-P_j-\bu)\theta(P-P_j+\bu)}{\theta^2(\bu) E(P_j, Q) E(P, P_j) d\zeta_j} =
\nn\\
&&
=\sum_{i=1}^n \res_{Z=P_i} z(Z)\, H_{PQ}(Z) \cdot \sum_{j=1}^n \res_{Z=P_j} z(Z)\, H_{QP}(Z) \, dz(P) dz(Q) =
\nn\\
&&
=\left[ z(P)-z(Q)\right]^2 \frac{\theta(P-Q-\bu)\theta(P-Q+\bu)}{\theta^2(\bu) E^2(P,Q) }.
\nn
\eeqa\epf

Equation \eqref{F2} immeaditely follows from Proposition \ref{prop223} and eq. \eqref{pref2}.

\setcounter{theorem}{0}
\subsection{Algebro-geometric solutions to the $n$-wave hierarchy and their tau-functions. Proof of eq. \eqref{F3}}\par

According to the original idea of S.P.Novikov \cite{Nov74} algebro-geometric (aka \emph{finite gap}) solutions to integrable systems of PDEs are obtained by considering stationary points of a linear combination of the commuting flows. Here we will be dealing with the $n$-wave system of nonlinear evolution PDEs represented in the form
\eqa
&&
\left[ L_{a,k}, L_{b,l}\right]=0
\nn\\
&&
L_{a,k}=\frac{\pal}{\pal t^a_k} -U_{a,k}({\bf t}; z), \quad a=1, \dots, n, \quad k\geq -1
\nn
\eeqa
where 
$$
U_{a,k}({\bf t}; z)=z^{k+1} E_a+\text{lower degree terms}
$$
is a $n\times n$ matrix-valued polynomial in $z$ of degree $k+1$ depending on the infinite number of independent variables $\bt=\left(t^a_k\right)$. The independent variable $z$ is often called \emph{spectral parameter}. The above equations hold true identically in $z$. Here the diagonal $n\times n$ matrix $E_a$ has only one nonzero entry
$$
\left( E_a\right)_{ij}=\delta_{ia}\delta_{a j}.
$$
For example, for $k=-1$
$$
U_{a,-1}=E_a
$$
and, for
$k=0$
$$
U_{a,0}=z\, E_a-\left[ E_a, Y\right], \quad a=1, \dots, n
$$
where 
the diagonal entries of the $n\times n$ matrix $Y=Y(\bt)$ vanish.
It turns out that the coefficients of the matrix polynomials $U_{a,k}(\bt; z)$ can be represented as polynomials in the entries of the matrix $Y(\bt)$ and its derivatives in the variables $t^1_0$, \dots, $t^n_0$. Thus the $n$-wave system can be considered as an infinite family of partial differential equations for the matrix $Y(\bt)$. See Appendix \ref{appa} below for the details about the structure of the $n$-wave hierarchy.

The following statement is crucial for computing tau-functions of solutions to the $n$-wave hierarchy.

\begin{prop} \label{mainprop} 1) For any solution to the $n$-wave hierarchy there exists a unique $n$-tuple of matrix-valued series 
$$
M_b(\bt,z)=E_b+\sum_{k\geq 1} \frac{B_{b,k}(\bt)}{z^k}, \quad b=1, \dots, n
$$
satisfying
$$
\left[L_{a,k}, M_b\right]=0\quad\Leftrightarrow\quad \frac{\pal M_b(\bt,z)}{\pal t^a_k}=\left[ U_{a,k}(\bt,z), M_b(\bt,z)\right]\quad \forall\, a=1,\dots, n, \quad k\geq -1
$$
and also
$$
M_a(\bt,z)M_b(\bt,z)=\delta_{ab}M_a(\bt,z), \quad M_1(\bt,z)+\dots+M_n(\bt,z)={\bf 1}.
$$

2) The (principle) tau-function $\tau(\bt)$ of this solution is determined from the following generating series in independent variables $z_1$, $z_2$ for its second logarithmic derivatives
\beq\label{maintau}
\sum_{p,\, q=0}^\infty  \frac1{z_1^{p+2}} \frac1{z_2^{q+2}} \frac{\pal^2\log\tau(\bt)}{\pal t^b_q \pal t^a_p}=\tr\frac{ M_a(\bt,z_1)M_b(\bt,z_2)}{(z_1-z_2)^2}-\frac{\delta_{ab}}{(z_1-z_2)^2}
\eeq
for any $a, \, b=1, \dots, n$.

3) The logarithmic derivatives of higher orders $N\geq 3$ of the same tau-function can be determined by the following generating series
\beq\label{maintau1}
\sum_{k_1, \dots, k_N\geq 0} \frac{\pal^N \log \tau(\bt)}{\pal t_{k_1}^{a_1}\dots \pal t_{k_N}^{a_N}}\frac1{z_1^{k_1+2}\dots z_N^{k_N+2}}=-\frac1{N}\sum_{s\in S_N}\frac{\tr\,\left[  M_{a_{s_1}}\left(\bt, z_{s_1}\right)\dots M_{a_{s_N}}\left(\bt, z_{s_N}\right)\right]}{\left(z_{s_1}-z_{s_2}\right)\dots \left( z_{s_{N-1}} -z_{s_N}\right) \left( z_{s_N}-z_{s_1}\right)}
\eeq
\end{prop}

Clearly the tau-function is determined by \eqref{maintau} uniquely up to
$$
\tau(\bt)\mapsto e^{\alpha+\sum \beta_{a,k}t^a_k} \tau(\bt)
$$
for some constants $\alpha$ and $\beta_{a,k}$. 

The above proposition about construction of tau-functions of solutions to the $n$-wave integrable hierarchy is an extension to this case of the approach of \cite{BDY1}--\cite{BDY3} based on the theory of the so-called matrix resolvents. For the proofs see Appendix \ref{appa} below.

\begin{remark} The construction of tau-function given in the Proposition differs from the original definition of \cite{Kac}, \cite{Di}, \cite{Wi} (we recall this definition in Appendix B below). One can proof equivalence of the two definitions following the scheme of \cite{BDY1}. We will not do it here as it is not needed for the proofs of the results of this paper.
\end{remark}

We will now apply the Proposition to the finite-gap solutions of the $n$-wave hierarchy. According to the Novikov's recipe mentioned above for arbitrary choice of constants $c_{a,k}$, $a=1, \dots, n$, $-1\leq k \leq N$ for any integer $N\geq 0$, we obtain a family of algebro-geometric solutions $Y(\bt)$ satisfying
$$
\sum_{a=1}^n\sum_{k=-1}^N c_{a,k} \frac{\pal Y(\bt)}{\pal t^a_k} =0.
$$
For any such solution define a matrix polynomial
$$
U(\bt,z)=\sum_{a=1}^n\sum_{k=-1}^N c_{a,k}U_{a,k}(\bt,z).
$$
Recall that the coefficients of the matrix polynomial belong to the space ${\mathcal Y}$.
From commutativity \eqref{nwoper} of the flows it readily follows that the matrix $U=U(\bt,z)$ satisfies
\beq\label{nov1}
\frac{\pal U}{\pal t^b_l}=\left[ U_{b,l},U\right]\quad \mbox{for any}\quad b=1, \dots, n, \quad l\geq -1.
\eeq
Therefore the characteristic polynomial
$$
R(z,w)=\det (w\cdot {\bf 1}-U(\bt,z))
$$
does not depend on $\bt$. Its coefficients can be considered as first integrals of the differential equations \eqref{nov1}. For a given $N$, assuming the coefficients $c_{1,N}$, \dots, $c_{n,N}$ to be pairwise distinct the matrix polynomial $U(\bt,z)$ for any $\bt$ belongs to the family ${\mathcal W}$ of polynomials of the form \eqref{spacenw} with $m=N+1$. Therefore the spectral curves $C=\{ R(z,w)=0\}$ is of the form \eqref{nwcurva}. 

Our nearest goal is to construct an algebro-geometric solution to the $n$-wave system such that the matrix polynomial $U(\bt,z)$ satisfies the initial condition
\beq\label{nov2}
W(z)=U(0,z)
\eeq
so $C$ is the spectral curve of the matrix polynomial $W(z)$.
The construction is rather standard for the theory of integrable systems. Namely, for the matrix polynomial $W(z)$ we have constructed the spectral curve $C$ with a nonspecial divisor $D_0$.  Starting from these data one can construct an algebro-geometric solution of the $n$-wave system following the I.M.Krichever's scheme \cite{Krich77}. We will use a  \emph{vector-valued Baker--Akhiezer function} $\bpsi(\bt,P)=\left( \psi^1(\bt, P), \dots, \psi^n(\bt,P)\right)^T$ meromorphic on the spectral curve $C\setminus\left(P_1\cup\dots\cup P_n\right)$ with poles at the points of the divisor $D_0$ and having essential singularities at $P=P_1$,\dots, $P=P_n$ of the form
$$
\psi^i(\bt,P)=\left(\delta_{ij}+{\mathcal O}\left(\frac1{z}\right)\right) e^{\phi_j(\bt,z)}, \quad P\to P_j,\quad i,\, j=1, \dots, n.
$$
It is a standard fact of the theory of Baker--Akhiezer functions that $\bpsi(\bt,P)$ exists for sufficiently small $|\bt|$ and is unique. It also exists the unique dual Baker--Akhiezer function $\bpsi^\dagger(\bt,P)=\left(\psi^\dagger_1(\bt,P),\dots, \psi^\dagger_n(\bt,P)\right)$ with poles at the divisor $D_0^\dagger$ and essential singularities of the form
$$
\psi^\dagger_i(\bt,P)=\left(\delta_{ij}+{\mathcal O}\left(\frac1{z}\right)\right) e^{-\phi_j(\bt,z)}, \quad P\to P_j,\quad i,\, j=1, \dots, n.
$$
The needed algebro-geometric solution to the $n$-wave system is uniquely specified by the condition that its wave function is expressed in terms of $\bpsi(\bt,P)$. We will now obtain an explicit expression of this solution in terms of theta-functions of the spectral curve.

\begin{prop} \label{prop251} Let $W(z)$ be a matrix polynomial of the form \eqref{spacenw} with a nonsingular spectral curve $C$. Let $D_0$ be the divisor of poles of the eigenvector of $W(z)$ normalized by the condition \eqref{nwnorm} and denote $(\bu_0, \bla^0)\in J(C, P_1,\dots, P_n)$ the corresponding point \eqref{corp} of the generalized Jacobian. Introduce Abelian differentials
\eqa\label{difomt2}
&&
\Omega^i_j(\bt,P)=\frac{\lambda_i(\bt)}{\lambda_j(\bt)}\frac{\theta(P-P_i-\bu(\bt))\theta(P-P_j+\bu(\bt))}{\theta^2(\bu(\bt))E(P_i,P)E(P,P_j)\sqrt{d\zeta_i}\sqrt{d\zeta_j}}
\\
&&
\nn\\
&&
\bu(\bt)=\bu_0-\sum t^a_k\bv^{(a,k)},\quad\lambda_i(\bt)=\exp\left\{ \sum t^a_k \dashint_{P_i}^{P_0}\Omega_a^{(k)}\right\}\lambda_i^0, \quad i, \, j=1, \dots, n.
\nn
\eeqa
Here and below
$$
\zeta_a=\zeta_a(Q)=\frac1{z(Q)}, \quad Q\in C, \quad Q\to P_a, \quad a=1, \dots, n
$$
is a natural local parameter near $P_a$. The principal values of the integrals are defined by the following limits
\beq\label{v.p.}
\dashint_{P_i}^{P_0} \Omega^{(k)}_a =\lim_{Q\to P_i}\left( \int_Q^{P_0}\Omega^{(k)}_a +z^{k+1}(Q)\right).
\eeq
Define matrix-valued power series in $z^{-1}$ by expanding the differentials at infinity
\beq\label{pita}
M_a(\bt,z)=\left(\frac{\Omega^i_j(\bt,P)}{dz}\right)_{1\leq i, \, j\leq n},\quad P=(z,w_a(z))\to P_a, \quad a=1, \dots, n
\eeq
and put
\beq\label{uak}
U_{a,k}(\bt,z)=\left(z^{k+1} M_a(\bt,z)\right)_+.
\eeq
This collection of matrix polynomials is an algebro-geometric solution to the $n$-wave system with the corresponding matrix $U(\bt,z)$ satisfying \eqref{nov1} given by
\beq\label{utz}
U(\bt,z)=w_1(z)M_1(\bt,z)+\dots+w_n(z)M_n(\bt,z).
\eeq
In this formula $w_a(z)$ is the Laurent expansion of the algebraic function $w(z)$ near $P=P_a$, $a=1$,\dots, $n$.
This matrix polynomial $U(\bt,z)$ satisfies the initial condition \eqref{nov2}.
\end{prop}

\pf Let $\Omega(P)$ be the differential \eqref{difom1} on $C$ constructed above. We first prove that the differentials \eqref{difomt2} coincide with
\beq\label{difomt1}
\Omega^i_j(\bt,P)=\psi^i(\bt,P)\Omega(P)\psi^\dagger_j(\bt,P),\quad i,\, j=1,\dots, n
\eeq
To this end we use the following expressions of the Baker--Akhiezer functions $\bpsi(\bt,P)$ and $\bpsi^\dagger(\bt,P)$
\beq\label{psii}
\psi^i(\bt,P)=\exp\left(\sum t^a_k ~ \dashint_{P_i}^P  \Omega^{(k)}_a\right)
\frac{
\lambda_i^0 
\frac{
\theta(P-P_i+\sum t^a_k \bv^{(a,k)} -\bu_0)
}
{
\theta(\sum t^a_k \bv^{(a,k)} -\bu_0) E(P,P_i) \left( d\zeta_i\right)^{1/2}}}
{\sum_{b=1}^n \lambda_b^0 \frac{\theta(P-P_b-\bu_0)}{\theta(\bu_0)E(P,P_b) \left( d\zeta_b\right)^{1/2}}},\quad i=1, \dots, n
\eeq
and
\beq\label{psidagi}
\psi^\dagger_i(\bt,P)=\exp\left(-\sum t^a_k ~ \dashint_{P_i}^P  \Omega^{(k)}_a\right)
\frac{
\frac1{\lambda_i^0} 
\frac{
\theta(P-P_i-\sum t^a_k \bv^{(a,k)} +\bu_0)
}
{
\theta(\sum t^a_k \bv^{(a,k)} -\bu_0) E(P_i,P) \left( d\zeta_i\right)^{1/2}}}
{\sum_{b=1}^n \frac1{\lambda_b^0} \frac{\theta(P-P_b+\bu_0)}{\theta(\bu_0)E(P_b,P) \left( d\zeta_b\right)^{1/2}}},\quad i=1, \dots, n.
\eeq
Here $\Omega^{(k)}_a$ is
the normalised second kind differential  on $C$ with a unique pole at $P_a$ of order $k+2$
\eqa\label{nwomega}
&&
\Omega^{(k)}_a(P)=dz^{k+1}+\mbox{regular terms}, \quad P\to P_a
\nn\\
&&
\oint_{a_i}\Omega^{(k)}_a=0, \quad V^{(a,k)}_i=\oint_{b_i}\Omega^{(k)}_a, \quad i=1, \dots, g
\eeqa
for a chosen canonical basis $a_i$, $b_j$ in $H_1(C, \mathbb Z)$. Recall \cite{Fay} that the $b$-periods $V_i^{(a,k)}$ of the differentials $\Omega^{(k)}_a$ coincide with the coefficients of expansions \eqref{nwholo} of holomorphic differentials $\omega_i(P)$ at $P\to P_a$.
Observe that for $\bt=0$ the functions $\bpsi(0,P)$ and $\bpsi^\dagger(0,P)$ coincide with the right and left eigenvectors of $W(z)$.

The derivation of the representations \eqref{psii}, \eqref{psidagi} is standard for the theory of Baker--Akhiezer functions: we check that \eqref{psii}, \eqref{psidagi} are well-defined meromorphic functions on $C\setminus\left(P_1\cup\dots\cup P_n\right)$ with essential singularities at infinity of the needed form having poles at rhe points of the divisors $D_0$ and $D^\dagger_0$ respectively (for the claim about the location of poles use Proposition \ref{prop221} and Remark \ref{rem225}). With the help of these expressions it is easy to verify validity of eq. \eqref{difomt1}. 

Define now a matrix-valued function $\Psi(\bt,z)$ such that its
$i$-th row is given by the expansions of $\psi^i(\bt,P)$ at the infinite points $P_1$, \dots, $P_n$. We will prove that $\Psi(\bt,z)$ is the wave function of the solution \eqref{uak} to the $n$-wave system. To this end we will first verify that the definition \eqref{pita} of the matrices $M_a(z)$ can be rewritten in the form \eqref{mpsi}. 

Introduce another matrix-valued function $\Psi^\dagger(\bt,z)$ in a similar way: its $i$-th column is given by the expansions of $\psi^\dagger_i(\bt,P)$ at the infinite points $P_1$, \dots, $P_n$. Let us prove that this matrix is inverse to $\Psi(\bt,z)$ up to multiplication on the left by a nondegenerate diagonal matrix.

\begin{lemma} \label{lem252}Define
\beq\label{omhat}
\hat\rho(z)=\diag\left( \rho_1(z), \dots, \rho_n(z)\right), \quad \rho_a(z)=\left(\frac{\Omega(P)}{dz}\right)_{P=(z,w_a(z))}=1+{\mathcal O}\left(\frac1{z}\right), \quad a=1, \dots, n.
\eeq
Then
\beq\label{inverse}
\Psi(\bt,z)\hat\rho(z)\Psi^\dagger(\bt,z)={\bf 1}.
\eeq
\end{lemma}

\pf The differentials \eqref{difomt2} are holomorphic on $C\setminus\left(P_1\cup\dots\cup P_n\right)$. At infinity they behave in the same way as the differentials $\Omega^i_j(P)$ (see Proposition \ref{prop217} above). For an arbitrary complex number $z$ away from the ramification points of $C$ consider the sum
$$
\Omega^i_j(\bt, (z,w_1(z))+\dots+\Omega^i_j(\bt,(z,w_n(z))=\Psi^i_1(\bt,z)\rho_1(z){\Psi^\dagger}^1_j(\bt,z)dz+\dots+\Psi^i_n(\bt,z)\rho_n(z){\Psi^\dagger}^n_j(\bt,z)dz.
$$
This is a well-defined differential on ${\bf P}^1$. It can have poles only at $z=\infty$, namely, a simple pole for $i\neq j$ and a double pole $\sim dz$ for $i=j$. Therefore the above sum is equal to $\delta^i_j dz$. \epf

\begin{cor} The matrix series $M_a(\bt,z)$ coincide with
\beq\label{mta}
M_a(\bt,z)=\Psi(\bt,z)E_a \Psi^{-1}(\bt,z), \quad a=1, \dots n.
\eeq
\end{cor}

\begin{lemma} 1) The matrix $\Psi(\bt,z)$ satisfies
\beq\label{lem1}
\frac{\pal }{\pal t^a_k}\Psi(\bt,z)=U_{a,k}(\bt,z) \Psi(\bt,z) \quad \forall~ a, \, b=1, \dots, n, \quad k\geq -1
\eeq
where $U_{a,k}(\bt,z)$ are given by \eqref{uak}.

\noindent 2) The matrix polynomials $U_{a,k}(\bt,z)$ satisfy eqs. \eqref{nwoper} of the $n$-wave hierarchy.

\noindent 3) The matrix series $M_b(\bt,z)$ satisfy
\beq\label{lem2}
\frac{\pal M_b(z)}{\pal t^a_k}=\left[ U_{a,k}(z),M_b(z)\right] .
\eeq
\end{lemma}

\pf Let $\Psi({\bf t},z)=A({\bf t}, z)e^{\phi({\bf t},z)}$ with $A(\bt,z)={\bf 1}+{\mathcal O}\left(\frac1{z}\right)$ and denote
$$
\tilde U_{a,k}(\bt,z)=\left(z^{k+1} M_a(\bt,z)\right)_-.
$$
It is a power series in $z^{-1}$. We have
$$
\frac{\pal}{\pal t^a_k} \Psi(\bt,z)-U_{a,k}(\bt,z)\Psi(\bt,z)=\left(\frac{\pal A(\bt,z)}{\pal t^a_k}\cdot A^{-1}(\bt,z)+\tilde U_{a,k}(\bt,z)\right)\Psi(\bt,z).
$$
As the expression in the parenthesis contains only negative powers of $z$, the right hand side is a Baker--Akhiezer function on the curve $C$ with the same divisor of poles and with expansion at infinity of the form
$$
\left(\frac{\pal A(\bt,z)}{\pal t^a_k}\cdot A^{-1}(\bt,z)+\tilde U_{a,k}(\bt,z)\right)\Psi(\bt,z)={\mathcal O}\left(\frac1{z}\right) e^{\phi(\bt,z)}.
$$
Hence this Baker--Akhiezer function identically vanishes. This proves the first part of Lemma.

The equations \eqref{nwoper} readily follow from the compatibility
$$
\frac{\pal }{\pal t^a_k}\frac{\pal }{\pal t^b_l}\Psi(\bt,z)=\frac{\pal }{\pal t^b_l}\frac{\pal }{\pal t^a_k}\Psi(\bt,z).
$$
Finally the eq. \eqref{lem2} follows from \eqref{mta} and \eqref{lem1}. \epf

In a similar way one can verify validity of eq. \eqref{nov1} for the matrix $U(\bt,z)$ defined by\eqref{utz}. It remains to prove that this matrix is polynomial in $z$ satisfying the initial condition \eqref{nov2}. To this end we consider the differentials
$$
U^i_j(\bt,z)dz=w_1(z)\Omega^i_j(\bt, (z,w_1(z))+\dots+w_n(z)\Omega^i_j(\bt,(z,w_n(z)).
$$
Like in the proof of Lemma \ref{lem252} this is a differential on ${\bf P}^1$ with poles only at infinity. Hence it must be a polynomial. Since $\Omega^i_j(0,P)=\Omega^i_j(P)$ we have $U(0,z)=W(z)$. The Proposition is proved. \epf

We are now ready to compute the tau-function of the algebro-geometric solution. Define numbers $q_{a,k;b,l}$ as coefficients of expansions of the second kind differentials $\Omega_a^{(k)}(P)$ at $P\to P_l$
\beq\label{qab}
\Omega_a^{(k)}(P)=\delta_{ab} d\left( z^{k+1}\right)+\sum_{l\geq 0}\frac{q_{a,k; b,l}}{z^{l+2}}dz, \quad P\to P_b.
\eeq
Alternatively these coefficients can be recovered from the expansions of the normalized second kind bi-differential \cite{Fay}
\eqa\label{bidif}
&&
\omega(P,Q)=d_Pd_Q\log E(P,Q)=\left[ \frac{\delta_{ab}}{(z_1-z_2)^2}+\sum_{k,\,l\geq 0}\frac{q_{a,k; b,l}}{z_1^{k+2}z_2^{l+2}}\right] dz_1 dz_2
\\
&&
z_1=z(P), ~z_2=z(Q),\quad P\to P_a, ~Q\to P_b.
\nn
\eeqa

\begin{prop} Tau-function of the algebro-geometric solution constructed in Proposition \ref{prop251} is equal to
\beq\label{taumain}
\tau(\bt)=e^{\frac12 \sum q_{a,k; b,l} t^a_k t^b_l}\theta\left(\sum t^a_k \bv^{(a,k)}-\bu_0\right)
\eeq
up to multiplication by exponential of a linear function. Here
\beq\label{vvec}
\bv^{(a,k)}=\left(V_1^{(a,k)},\dots, V_g^{(a,k)}\right).
\eeq
\end{prop}

\pf We have to compute the generating function \eqref{cor2} of the second logarithmic derivatives of the tau-function
\beq\label{step0}
\sum\frac{\pal^2\log\tau(\bt)}{\pal t^a_k \pal t^b_l}\frac{dz_1}{z_1^{k+2}}\frac{dz_2}{z_2^{l+2}}=\frac{\tr [M_a(\bt,z_1)M_b(\bt,z_2)]}{(z_1-z_2)^2}dz_1 dz_2 -\frac{\delta_{ab}}{(z_1-z_2)^2}dz_1 dz_2
\eeq
where $M_1(\bt,z)$, $M_b(\bt,z)$ are solutions to the equations \eqref{lem2} in the class of matrix-valued power series in $z^{-1}$ uniquely specified by the conditions
\eqref{mta}. Using the representation \eqref{pita} we can rewrite the previous equation in the form
\eqa
&&
\sum\frac{\pal^2\log\tau(\bt)}{\pal t^a_k \pal t^b_l}\frac{dz_1}{z_1^{k+2}}\frac{dz_2}{z_2^{l+2}}=\frac{\tr [\bom(\bt,P)\bom(\bt,Q)]}{(z(P)-z(Q))^2}-\frac{\delta_{ab}}{(z_1-z_2)^2}dz_1 dz_2
\nn\\
&&
\nn\\
&&
 z_1=z(P),~ z_2=z(Q), \quad P\to P_a, ~ Q\to P_b
\nn
\eeqa
where the matrix entries of $\bom(\bt,P)$ are equal to \eqref{difomt2}. According to Proposition \ref{prop226} the right hand side can be rewritten in the form
\beq\label{step}
\frac{\tr [\bom(\bt,P)\bom(\bt,Q)]}{(z(P)-z(Q))^2}-\frac{\delta_{ab}}{(z_1-z_2)^2}dz_1 dz_2=\frac{\theta(P-Q-\bu(\bt))\theta(P-Q+\bu(\bt))}{\theta^2(\bu(\bt))E^2(P,Q)}-\frac{\delta_{ab}}{(z_1-z_2)^2}dz_1 dz_2.
\eeq
We will now use the following important identity \cite{Fay}
$$
\frac{\theta(P-Q-\bu(\bt))\theta(P-Q+\bu)}{\theta^2(\bu)E^2(P,Q)}=\sum_{i,\,j=1}^n \frac{\pal^2 \log\theta(\bu)}{\pal u_i \pal u_j} \omega_i(P)\omega_j(Q)+\omega(P,Q)
$$
where $\omega(P,Q)$ is the normalized bi-differential \eqref{bidif}. Using this identity we can expand the right hand side of eq. \eqref{step} at $P\to P_a$, $Q\to P_b$ 
\eqa
&&
\frac{\theta(P-Q-\bu(\bt))\theta(P-Q+\bu(\bt))}{\theta^2(\bu(\bt))E^2(P,Q)}-\frac{\delta_{ab}}{(z_1-z_2)^2}dz_1 dz_2
\nn\\
&&
=\sum_{i,\,j=1}^n \frac{\pal^2 \log\theta(\bu(\bt))}{\pal u_i \pal u_j} \sum_{k\geq 0} \frac{V_i^{(a,k)}}{z_1^{k+2}}\sum_{l\geq 0} \frac{V_j^{(b,l)}}{z_2^{l+2}}dz_1dz_2+\sum_{k,\,l\geq 0}\frac{q_{a,k; b,l}}{z_1^{k+2}z_2^{k+2}} dz_1 dz_2
\nn\\
&&
=\sum_{k,\,l\geq 0} \frac{\pal^2\log\theta(\bu(\bt))}{\pal t^a_k \pal t^b_l}\frac{dz_1}{z_1^{k+2}}\frac{dz_2}{z_2^{l+2}}+\sum_{k,\,l\geq 0}\frac{q_{a,k; b,l}}{z_1^{k+2}z_2^{k+2}} dz_1 dz_2.
\nn
\eeqa
Comparing this expansion with \eqref{step0} we arrive at the proof of the Proposition. \epf

\begin{remark} An expression similar to \eqref{taumain} is well known in the theory of KP equation and its reductions \cite{Segal}, \cite{EH}, \cite{Naka}. We emphasize that here our main task was to prove that eq. \eqref{taumain} is in agreement with the construction of the tau-function given in terms of Proposition \ref{mainprop}.
\end{remark}

Let us now proceed to the proof of eq. \eqref{F3}. The expression \eqref{corN} for the $N$-th order logarithmic derivatives of the tau-function will be applied to the tau-function \eqref{taumain} of an algebro-geometric solution. Due to the previous Proposition the tau-function in the left hand side of eq. \eqref{corN} for $N\geq 3$ can be replaced with the theta-function
\eqa
&&
\sum_{k_1, \dots, k_N\geq 0} \frac{\pal^N \log \tau(\bt)}{\pal t_{k_1}^{a_1}\dots \pal t_{k_N}^{a_N}}\frac{dz_1\dots dz_N}{z_1^{k_1+2}\dots z_N^{k_N+2}}=\sum_{k_1, \dots, k_N\geq 0} \frac{\pal^N \log \theta(\bu(\bt))}{\pal t_{k_1}^{a_1}\dots \pal t_{k_N}^{a_N}}\frac{dz_1\dots dz_N}{z_1^{k_1+2}\dots z_N^{k_N+2}}=
\nn\\
&&
=(-1)^N\sum_{k_1, \dots, k_N\geq 0}\sum_{i_1,\dots,i_N=1}^gV_{i_1}^{(a_1,k_1)}\dots V_{i_N}^{(a_N,k_N)} \frac{\pal^N \log \theta(\bu(\bt))}{\pal u_{i_1}\dots \pal u_{i_N}}\frac{dz_1\dots dz_N}{z_1^{k_1+2}\dots z_N^{k_N+2}}=
\nn\\
&&
=(-1)^N\sum_{i_1,\dots,i_N=1}^g\frac{\pal^N \log \theta(\bu(\bt))}{\pal u_{i_1}\dots \pal u_{i_N}}\omega_{i_1}(Q_1)\dots \omega_{i_N}(Q_N), \quad Q_1\to P_{a_1},\dots, Q_N\to P_{a_N}
\nn
\eeqa
where the last multi-differential is considered as its expansion in negative powers of
$$
z_1=z(Q_1), \dots, z_N=z(Q_N).
$$

Let us now consider the right hand side of eq. \eqref{corN} multiplying it, like above by $dz_1\dots dz_N$
\eqa
&&
-\frac1{N}\sum_{s\in S_N}\frac{\tr\,\left[  M_{a_{s_1}}\left(\bt, z_{s_1}\right)\dots M_{a_{s_N}}\left( \bt, z_{s_N}\right)\right]}{\left(z_{s_1}-z_{s_2}\right)\dots \left( z_{s_{N-1}} -z_{s_N}\right) \left( z_{s_N}-z_{s_1}\right)}dz_1\dots dz_N=
\nn\\
&&
=-\frac1{N}\sum_{s\in S_N}\frac{\tr\,\left[ \bom\left(\bt, Q_{s_1}\right)\dots\bom\left(\bt, Q_{s_N}\right)\right]}{\left(z(Q_{s_1})-z(Q_{s_2})\right)\dots \left(z(Q_{s_{N-1}})-z(Q_{s_N})\right)\left(z(Q_{s_N})-z(Q_{s_1})\right) }
\nn
\eeqa
where $Q_1$, \dots, $Q_N$ are arbitrary points of $C$ such that $Q_1\to P_{a_1}$, \dots $Q_N\to P_{a_N}$. So eq. \eqref{corN} implies that the two $N$-differentials coincide when the points $Q_1$, \dots, $Q_N$ go to infinity in all possible ways. Therefore these $N$-differentials coincide
\eqa\label{F3new}
&&
\sum_{i_1,\dots,i_N=1}^g\frac{\pal^N \log \theta(\bu(\bt))}{\pal u_{i_1}\dots \pal u_{i_N}}\omega_{i_1}(Q_1)\dots \omega_{i_N}(Q_N)=
\\
&&
=\frac{(-1)^{N-1}}{N}\sum_{s\in S_N}\frac{\tr\,\left[ \bom\left(\bt, Q_{s_1}\right)\dots\bom\left(\bt, Q_{s_N}\right)\right]}{\left(z(Q_{s_1})-z(Q_{s_2})\right)\dots \left(z(Q_{s_{N-1}})-z(Q_{s_N})\right)\left(z(Q_{s_N})-z(Q_{s_1})\right) }.
\nn
\eeqa
To complete the derivation of eq. \eqref{F3} we set $\bt=0$ where
$$
\bu(0)=\bu_0, \quad \bom(0,P)=\Phi(P)\frac{dz}{R_w(z,w)}.
$$
The Main Lemma and Main Theorem are proved. \epf

\setcounter{theorem}{0}
\subsection{Proof of Corollary \ref{cor03}}\par

Let $C$ be a compact Riemann surface of genus $g>0$ and $n$, $m$ a pair of positive integers.

\begin{prop} For sufficienly large $n$, $m$ and an arbitrary collection of $n$ pairwise distinct points $P_1$,\dots, $P_n$ there exist two rational functions $z$, $w$ on $C$ such that

(i) the function $z$ has simple poles at $P_1$,\dots, $P_n$

(ii) the function $w$ has poles of order $m$ at the same points and $w\neq P(z)$ for any polynomial $P$.
\end{prop} 

This is an easy consequence of Riemann--Roch theorem.

\begin{cor} An arbitrary compact Riemann surface of genus $g>0$ can be represented as the spectral curve of a matrix $W(z)$ of the form \eqref{mapol} for sufficiently large $n$ and $m$.
\end{cor}

According to the Corollary we can rewrite eq. \eqref{F3} in the form
\eqa\label{F3new1}
&&
\sum_{i_1,\dots,i_N=1}^g\frac{\pal^N \log \theta(\bu)}{\pal u_{i_1}\dots \pal u_{i_N}}\omega_{i_1}(Q_1)\dots \omega_{i_N}(Q_N)=
\\
&&
=\frac{(-1)^{N-1}}{N}\sum_{s\in S_N}\frac{\tr\,\left[ \bom\left( Q_{s_1}\right)\dots\bom\left( Q_{s_N}\right)\right]}{\left(z(Q_{s_1})-z(Q_{s_2})\right)\dots \left(z(Q_{s_{N-1}})-z(Q_{s_N})\right)\left(z(Q_{s_N})-z(Q_{s_1})\right) }.
\nn
\eeqa
(cf. \eqref{F3new}). Here $\bu\in J(C)\setminus (\theta)$ is the point of the Jacobian corresponding to the matrix $W(z)$, the matrix-valued differential $\bom(Q)$ equals
$$
\bom(Q)=\Phi(Q)\frac{dz}{R_w(z,w)}.
$$
Using the representation \eqref{matom} of the matrix entries of this differential we can rewrite the numerator, for an arbitrary permutation $s\in S_N$ as follows
$$
\tr\,\left[ \bom\left( Q_{s_1}\right)\dots\bom\left( Q_{s_N}\right)\right]=
$$
$$
=\sum_{i_1,\dots,i_N=1}^n
\frac{\theta(Q_{s_1}-P_{i_1}-\bu)\theta(Q_{s_1}-P_{i_2}+\bu)}{\theta^2(\bu) E(P_{i_1},Q_{s_1}) E(Q_{s_1},P_{i_2})\sqrt{d\zeta_{i_1}}\sqrt{d\zeta_{i_2}}}\frac{\theta(Q_{s_2}-P_{i_2}-\bu)\theta(Q_{s_2}-P_{i_3}+\bu)}{\theta^2(\bu) E(P_{i_2},Q_{s_2}) E(Q_{s_2},P_{i_3})\sqrt{d\zeta_{i_2}}\sqrt{d\zeta_{i_3}}}\dots
$$
$$
\dots 
\frac{\theta(Q_{s_{N-1}}-P_{i_{N-1}}-\bu)\theta(Q_{s_{N-1}}-P_{i_N}+\bu)}{\theta^2(\bu) E(P_{i_{N-1}},Q_{s_{N-1}}) E(Q_{s_{N-1}},P_{i_N})\sqrt{d\zeta_{i_{N-1}}}\sqrt{d\zeta_{i_N}}}
\frac{\theta(Q_{s_N}-P_{i_N}-\bu)\theta(Q_{s_N}-P_{i_1}+\bu)}{\theta^2(\bu) E(P_{i_N},Q_{s_N}) E(Q_{s_N},P_{i_N})\sqrt{d\zeta_{i_N}}\sqrt{d\zeta_{i_1}}}
$$
$$
= 
 \sum_{i_2=1}^n\frac{\theta(Q_{s_2}-P_{i_2}-\bu)\theta(Q_{s_1}-P_{i_2}+\bu)}{\theta^2(\bu) E(P_{i_2}, Q_{s_2}) E(Q_{s_1}, P_{i_2})  \,d\zeta_{i_2}}
  \sum_{i_3=1}^n\frac{\theta(Q_{s_3}-P_{i_3}-\bu)\theta(Q_{s_2}-P_{i_3}+\bu)}{\theta^2(\bu) E(P_{i_3}, Q_{s_3}) E(Q_{s_2}, P_{i_3})  \,d\zeta_{i_3}}
 \dots
$$
$$
\dots\sum_{i_N=1}^n\frac{\theta(Q_{s_N}-P_{i_N}-\bu)\theta(Q_{s_{N-1}}-P_{i_N}+\bu)}{\theta^2(\bu) E(P_{i_N}, Q_{s_N}) E(Q_{s_{N-1}}, P_{i_N})  \,d\zeta_{i_N}}\sum_{i_1=1}^n\frac{\theta(Q_{s_1}-P_{i_1}-\bu)\theta(Q_{s_N}-P_{i_1}+\bu)}{\theta^2(\bu) E(P_{i_1}, Q_{s_1}) E(Q_{s_N}, P_{i_1})  \,d\zeta_{i_1}} 
$$
$$
=
\sum_{i_2=1}^n \res_{Z=P_{i_2}}z(Z)H_{Q_{s_2} Q_{s_1}}(Z)
\dots
\sum_{i_N=1}^n \res_{Z=P_{i_N}}z(Z)H_{Q_{s_N} Q_{s_{N-1}}}(Z)\sum_{i_1=1}^n \res_{Z=P_{i_1}}z(Z)H_{Q_{s_1} Q_{s_N}}(Z)=
$$
$$
=(z_{s_2}-z_{s_1})\dots (z_{s_N}-z_{s_{N-1}})(z_{s_1}-z_{s_N})\frac{\theta(Q_{s_2}-Q_{s_1}-\bu)\dots \theta(Q_{s_N}-Q_{s_{N-1}}-\bu)\theta(Q_{s_1}-Q_{s_N}-\bu)
}{\theta^N(\bu) E(Q_{s_1},Q_{s_2})\dots E(Q_{s_{N-1}},Q_{s_N}) E(Q_{s_N},Q_{s_1})
}.
$$
In the above computation we have used the differential $H_{PQ}$ defined by \eqref{hpq}. The computation of residues is analogous to the one in \eqref{sumid}. In the last line we use the short notation $z_s:=z(Q_s)$. The Corollary \ref{cor03} is proved. \epf

For $N=3$ the identity \eqref{id3} takes the following explicit form
\eqa\label{ex3}
&&
\theta^3(\bu) \sum_{i, \, j, \, k=1}^g \frac{\pal^3\log\theta(\bu)}{\pal u_i \pal u_j \pal u_k} \omega_i(Q_1) \omega_j(Q_2) \omega_k(Q_3)=
\\
&&
\frac{\theta(Q_1-Q_2-\bu)\theta(Q_2-Q_3-\bu)\theta(Q_3-Q_1-\bu)-\theta(Q_1-Q_2+\bu)\theta(Q_2-Q_3+\bu)\theta(Q_3-Q_1+\bu)}{E(Q_1,Q_2)E(Q_2, Q_3) E(Q_3, Q_1)}
\nn
\eeqa
and for $N=4$
\eqa\label{ex4}
&&
\theta^4(\bu) \sum_{i, \, j, \, k, \, l=1}^g \frac{\pal^4\log\theta(\bu)}{\pal u_i \pal u_j \pal u_k \pal u_l} \omega_i(Q_1) \omega_j(Q_2) \omega_k(Q_3)\omega_l(Q_4)=
\\
&&
=V_\bu(Q_1, Q_2, Q_3, Q_4)+V_\bu(Q_1, Q_3, Q_2, Q_4)+V_\bu(Q_1, Q_3, Q_4, Q_2)
\nn
\eeqa
where
$$
V_\bu(Q_1, Q_2, Q_3, Q_4)=-\frac{\theta(Q_1-Q_2+\bu)\dots\theta(Q_4-Q_1+\bu)+\theta(Q_1-Q_2-\bu)\dots\theta(Q_4-Q_1-\bu)}{E(Q_1,Q_2)E(Q_2, Q_3) E(Q_3, Q_4)E(Q_4, Q_1)}.
$$

\section{Examples}
\setcounter{equation}{0}
\setcounter{theorem}{0}

\subsection{Hyperelliptic case}
\setcounter{theorem}{0}

Consider a $2\times 2$ matrix polynomial of the form
\beq\label{nlsw}
W(z)=\left(\begin{array}{cc} a(z) & b(z)\\ c(z) & -a(z)\end{array}\right), \quad \deg a(z)=g+1, \quad \deg b(z)=\deg c(z)=g,
\eeq
the polynomial $a(z)$ is monic.
The spectral curve
\beq\label{nlscurva}
w^2=Q(z), \quad Q(z)=-\det W(z)=z^{2g+2}+q_1 z^{2g+1}+\dots+q_{2g+2}
\eeq
is hyperelliptic. It has two distinct points $P_\pm$ at infinity,
$$
w=\pm z^{g+1}+\dots, \quad (z,w)\to P_\pm.
$$
We have
\beq\label{primer2}
\Pi(z,w)=\frac12 \frac{w+W(z)}{w}
\eeq
so the basic idempotents of the matrix $W(z)$ take the form
$$
M_\pm(z) =\frac12 \pm \frac12 \frac{W(z)}{w(z)}, \quad w(z)=z^{g+1}\sqrt{1+\frac{q_1}{z}+\dots+\frac{q_{2g+2}}{z^{2g+2}}}.
$$
Thus eq. \eqref{F2} for the bi-differential takes the following form
$$
\frac{\theta\left( Q_1-Q_2-\bu_0\right)\theta\left( Q_1-Q_2+\bu_0\right)}{\theta^2(\bu_0) E(Q_1,Q_2)^2}=\frac{b(z_1)c(z_2)+2a(z_1)a(z_2)+b(z_2)c(z_1)+2 w_1 w_2}{4(z_1-z_2)^2w_1w_2} \,dz_1\,dz_2,
$$
$Q_1=(z_1,w_1)$, $Q_2=(z_2,w_2)$.

Since
$$
M_+(z)-M_-(z)=\frac{W(z)}{w(z)}
$$
and the time-derivatives satisfy
$$
\frac{\pal}{\pal t^+_k}+\frac{\pal}{\pal t^-_k}=0
$$
we introduce
$$
\frac{\pal}{\pal t_k}=\frac{\pal}{\pal t^+_k}-\frac{\pal}{\pal t^-_k}.
$$
So eq. \eqref{corN} for the tau-function of the spectral curve \eqref{nlscurva} reduces to
\beq\label{nlsglav}
\sum_{k_1, \dots, k_N}\frac{\frac{\pal^N\log\tau(0)}{\pal t_{k_1}\dots \pal t_{k_N}}}{z_1^{k_1+2}\dots z_N^{k_N+2}}=-\frac1{N}\frac1{w(z_1)\dots w(z_N)}\sum_{s\in S_N}\frac{\tr\left[ W(z_{s_1})\dots W(z_{s_N})\right]}{(z_{s_1}-z_{s_2})\dots (z_{s_{N-1}}-z_{s_N})(z_{s_N}-z_{s_1})}-\frac{2\delta_{N,2}}{(z_1-z_2)^2}
\eeq
for any $N\geq 2$. First few logarithmic derivatives $F_{i_1\dots i_N}:=\pal^N\log\tau(0)/\pal t_{k_1}\dots \pal t_{k_N}$ read
\eqa
&&
F_{00}=-b_1c_1,\quad
F_{01}=2 a_1 b_1 c_1 - b_2 c_1 - b_1 c_2,
\nn\\
&&
F_{11}=\frac12 (-8 a_1^2 b_1 c_1 + 4 a_2 b_1 c_1 + 6 a_1 b_2 c_1 - 
   2 b_3 c_1 + b_1^2 c_1^2 + 6 a_1 b_1 c_2 - 4 b_2 c_2 - 
   2 b_1 c_3)
   \nn\\
   &&
 F_{000}=2(b_1c_2-b_2c_1),\quad
 F_{001}=2 (a_1 b_2 c_1 - b_3 c_1 - a_1 b_1 c_2 + b_1 c_3)
 \nn\\
 &&
 F_{0000}=4 (2 a_2 b_1 c_1 - a_1 b_2 c_1 - b_3 c_1 - a_1 b_1 c_2 + 
   2 b_2 c_2 - b_1 c_3).
   \nn
\eeqa
We do not specify the genus: the above expressions are valid for any $g\geq 2$; for $g=1$ one has to set $a_3=b_3=c_3=0$.

The theta-function of the spectral curve is related to the tau-function by the equation\footnote{Like above the eq. \eqref{taunls} holds true up to multiplication by exponential of a linear function of $t_i$.}
\beq\label{taunls}
\tau(\bt)=e^{\frac12\sum_{i,\,j}q_{ij}t_i t_j}\theta\left(\sum_k t_k \bv^{(k)}-\bu_0\right), \quad \bt=(t_0, t_1,\dots )
\eeq
with suitable coefficients $q_{ij}$ (cf. eq. \eqref{bidif} above). The vectors $\bv^{(k)}=\left( V_1^{(k)},\dots, V_g^{(k)}\right)$ have the form 
\beq\label{hyperv}
V_i^{(k)}=\alpha_{i1}r_k+\alpha_{i2}r_{k-1}+\dots+\alpha_{ig}r_{k-g+1}, \quad i=1, \dots, g, \quad k\geq 0
\eeq
where the $g\times g$ matrix $\left(\alpha_{ij}\right)$ is the inverse, up to a factor $2\pi\sqrt{-1}$ to the matrix of $a$-periods of the following holomorphic differentials
\beq\label{hyperalpha}
\left(\alpha_{ij}\right)=2\pi \sqrt{-1}\left(\oint_{a_j} z^{g-i}\frac{dz}{2w}\right)^{-1}
\eeq
and the rational numbers $r_k$ come from the expansion
$$
\left(1+\frac{q_1}{z}+\dots+\frac{q_{2g+2}}{z^{2g+2}}\right)^{-1/2}=\sum_{k\geq 0 }\frac{r_k}{z^k};
$$

The point $\bu_0$ is given by
\beq\label{nlspoint}
\bu_0=\sum_{j=1}^{g+1}\left(\int_{P_+}^{Q_j}\omega_1,\dots, \int_{P_+}^{Q_j}\omega_g\right)-{\boldsymbol \varpi}
\eeq
\beq\label{hyperholo}
\omega_i=(\alpha_{i1}z^{g-1}+\dots+\alpha_{ig})\frac{dz}{2w}, \quad i=1, \dots, g
\eeq
and the half-period ${\boldsymbol \varpi}$ 
for a suitable choice of the basis of cycles (see details in \cite{Fay}) has the form
\beq\label{varpi}
{\boldsymbol \varpi}=\pi\sqrt{-1}(1,0,1,0,\dots)+\frac12\sum_{i=1}^g \left(B_{1i},B_{2i}, \dots, B_{gi}\right).
\eeq
Points of the divisor $D_0=Q_1+\dots+Q_{g+1}$ of poles of the normalized eigenvector of the matrix $W(z)$ have the form $Q_i=(z_i,w_i)$ where $z_1$, \dots, $z_{g+1}$ are roots of the equation
\beq\label{nlspoint1}
a(z)=\frac12(b(z)+c(z))
\eeq
and
\beq\label{nlspoint2}
w_i=\frac12(c(z_i)-b(z_i)),\quad i=1, \dots, g+1.
\eeq

The Corollary \ref{cor03} in this particular case takes the following form

\begin{cor} \label{cor411} Assume rationality of coefficients of the polynomials $a(z)$, $b(z)$, $c(z)$. Then for any $N\geq 3$ and an arbitrary choice of indices $k_1$, \dots, $k_N\geq 0$ one has
$$
\sum_{i_1, \dots, i_N=1}^g V_{i_1}^{(k_1)}\dots V_{i_N}^{(k_N)}\frac{\pal^N\log\theta(\bu_0)}{\pal u_{i_1}\dots \pal u_{i_N}}\in\mathbb Q.
$$
Here $\theta=\theta(\bu | B)$ is the theta-function \eqref{theta} of the hyperelliptic curve \eqref{nlscurva} and
the point $\bu_0$ is given by \eqref{nlspoint}, \eqref{nlspoint1}, \eqref{nlspoint2}.
 \end{cor}

\subsection{Three-sheet Riemann surfaces}
\setcounter{theorem}{0}

Let
\beq\label{3sheet}
W(z)=z^m B^0+z^{m-1} B^1+\dots +B^m,\quad B^k=\left( B_{ij}^k\right)_{1\leq i, j\leq 3}, \quad B^0=\diag (b^0_1,b^0_2, b^0_3)
\eeq
be a 3$\times$3 matrix polynomial satisfying $b^0_i\neq b^0_j$ for $i\neq j$ and $\tr\, W(z)=0$. Let
$$
w^3+p(z)w+q(z)=\det(w\cdot {\bf 1}-W(z))
$$
be the characteristic polynomial. The genus of the spectral curve $C$ is equal to $g=3m-2$. The spectral projectors of $W(z)$ are given by branches of the algebraic function
\beq\label{3pi}
\Pi(z,w)=\frac{W^2+w\,W +w^2 +p(z)}{3w^2+p(z)},\quad (z,w)\in C.
\eeq
The expession for the bidifferential \eqref{F2} takes the following form
\eqa\label{3bidi}
&&
\frac{\theta\left( Q_1-Q_2-\bu_0\right)\theta\left( Q_1-Q_2+\bu_0\right)}{\theta^2(\bu_0) E(Q_1,Q_2)^2}=
\\
&&
\frac{\tr\, \left[W_1^2 W_2^2+(w_1 W_2+w_2 W_1)W_1 W_2+w_1 w_2 W_1 W_2\right]-2(p_1 p_2 +w_1^2 p_2 + w_2^2 p_1 +3 w_1^2 w_2^2) }{(3w_1^2+p_1)(3w_2^2+p_2)(z_1-z_2)^2}dz_1dz_2
\nn
\eeqa
where $Q_i=(z_i ,w_i)$ and we use short notations
$$
W_i=W(z_i), \quad p_i=p(z_i),\quad i=1,\, 2.
$$
The first few logarithmic derivatives of the tau-function \eqref{idprin} read as follows
\eqa
&&
F^{ij}_{00}[W]=\frac{b^1_{ij}b^1_{ji}}{(b^0_i-b^0_j)^2},\quad i\neq j
\nn\\
&&
F^{ij}_{10}[W]=\frac{b^2_{ij}b^1_{ji}+b^2_{ji}b^1_{ij}}{(b^0_i-b^0_j)^2}-2(b^1_{ii}-b^1_{jj})\frac{b^1_{ij}b^1_{ji}}{(b^0_i-b^0_j)^3}
,\quad i\neq j
\nn\\
&&
F^{123}_{000}[W] =\frac{b^1_{12}b^1_{23}b^1_{31}-b^1_{13}b^1_{32}b^1_{21}}{(b^0_1-b^0_2)(b^0_2-b^0_3)(b^0_3-b^0_1)}
\nn\\
&&
F^{iij}_{000}[W]=\frac{b^1_{ij}b^1_{jk}b^1_{ki}-b^1_{ik}b^1_{kj}b^1_{ji}}{(b^0_i-b^0_j)^2(b^0_k-b^0_i)}+\frac{b^2_{ij}b^1_{ji}-b^2_{ji}b^1_{ij}}{(b^0_i-b^0_j)^2},\quad i\neq j, \quad k\neq i, \, j
\nn
\eeqa
etc. One can also compute the derivatives of the above type for $i=j$ using the identities 
$$
\sum_{a=1}^n \frac{\pal}{\pal t^a_k}=0, \quad \forall~k\geq 0.
$$

\appendix

\section{Appendix. Tau-function of the $n$-wave integrable system}\label{appa}
\setcounter{equation}{0}
\setcounter{theorem}{0}


The (complexified) $n$-wave system \cite{soliton} is an infinite family of pairwise commuting systems of nonlinear PDEs for $n(n-1)$ functions $y_{ij}$, $i\neq j$ of infinite number of independent variables $t^a_k$, $a=1, \dots, n$, $k\geq 0$ called \emph{times}. We will often use an alternative notation for the variables $t^a_0=:x^a$, $a=1, \dots, n$ that will be called \emph{spatial variables}.

The equations of the $n$-wave hierarchy (also called AKNS-D hierarchy, see \cite{Di2}) are written as conditions of commutativity of linear differential operators
\eqa\label{nwoper}
&&
L_{a,k}=\frac{\pal}{\pal t^a_k} -U_{a,k}({\bf y}; z), \quad a=1, \dots, n, \quad k\geq 0
\nn\\
&&
\left[ L_{a,k}, L_{b,l}\right]=0
\eeqa
where $U_{a,k}({\bf y}; z)$ is a $n\times n$ matrix-valued polynomial in $z$ of degree $k+1$ depending polynomially on the functions $y_{ij}$ and their derivatives in $x^1$, \dots, $x^n$. For $k=0$ one has 
\beq\label{nwx}
U_{a,0}=z\, E_a-\left[ E_a, Y\right]
\eeq
where the diagonal $n\times n$ matrix $E_a$ has only one nonzero entry
$$
\left( E_a\right)_{ij}=\delta_{ia}\delta_{a j},
$$
the $n\times n$ matrix $Y$ has the form
$$
Y=\left(y_{ij}\right), \quad y_{ii}=0.
$$
The commutativity
\beq\label{nw00}
\left[ L_{a,0}, L_{b,0}\right]=0
\eeq
implies the system of constraints
\eqa\label{nw0}
&&
\sum_{k=1}^n \frac{\pal y_{ij}}{\pal x^k}=0
\nn\\
&&
\frac{\pal y_{ij}}{\pal x^k} =y_{ik}y_{kj}, \quad \mbox{the indices} \quad i, \, j, \, k \quad \mbox{are pairwise distinct}.
\eeqa
For $n=2$ the second part of the constraints is empty. 

In order to construct the matrix polynomials $U_{a,k}$ for $k>0$ we will use the following procedure \cite{Dub77} that can be considered as a generalization of the well-known AKNS construction developed for $n=2$ in the seminal paper \cite{AKNS}. 

Consider an arbitrary function\footnote{Here and below saying ``functions" we have in mind just formal power series in the independent variables.}  $Y(\bx)$ satisfying the system \eqref{nw0}. We are looking for solutions to the following system of linear differential equations 
\beq\label{adj}
\frac{\pal M}{\pal x^a}=\left[ U_{a,0},M\right] \quad\Leftrightarrow\quad \left[ L_{a,0},M\right]=0,\quad a=1, \dots,n.
\eeq
for a matrix-valued function $M=M(\bx,z)$ of the form
\beq\label{riadm}
M(\bx,z)=\sum_{k\geq 0}\frac{M_k(\bx)}{z^k}.
\eeq
Compatibility of this overdetermined system of differential equations follows from \eqref{nw00}.
Observe that the coefficients of the characteristic polynomial $\det\left( M(\bx,z)-w\cdot{\bf 1}\right)$ of the matrix $M(\bx,z)$ are first integrals of the system \eqref{adj}.

\begin{prop} \label{prop242} For an arbitrary solution $Y(\bx)$ to the system \eqref{nw00}, \eqref{nw0} there exist unique matrix series of the form
\beq\label{max}
M_a(\bx,z)=E_a+\sum_{k\geq 1} \frac{B_{a,k}(\bx)}{z^k}, \quad a=1, \dots, n
\eeq
satisfying \eqref{adj} as well as the following equations
\beq\label{idem}
M_a(\bx,z) M_b(\bx,z)=\delta_{ab}M_a(\bx,z), \quad M_1(\bx,z)+\dots +M_n(\bx,z)={\bf 1}.
\eeq
\end{prop}

\pf We begin with the recursion procedure for computing the coefficients of the expansion \eqref{riadm}. Clearly $M_0$ must be a constant diagonal matrix. 
Other coefficients can be determined by the following procedure.

\begin{lemma} For an arbitrary solution $Y=Y(\bx)$ to eqs. \eqref{nw0} and an arbitrary diagonal matrix $B=\diag (b_1,\dots, b_n)$ with pairwise distinct diagonal entries there exists a unique solution
\beq\label{adj1}
M=M_B(\bx,z)=B+\sum_{k\geq 1} \frac{M_{B,k}(\bx)}{z^k}
\eeq
to the system \eqref{adj} normalized by the condition
\beq\label{nwnorma}
\det\left( M_B(\bx, z)-w\cdot {\bf 1}\right)=\det\left( B-w\cdot {\bf 1}\right).
\eeq
\end{lemma}

\pf Split every coefficient into its diagonal and off-diagonal part
$$
M_{B,k}=D_k+C_k, \quad k\geq 1.
$$
Vanishing of the constant term in \eqref{adj} implies
$$
[E_a,M_{B,1}]=\left[B,[E_a,Y]\right]=\left[E_a,[B,Y]\right],\quad a=1, \dots, n.
$$
This system uniquely determines the off-diagonal part of the matrix $M_{B,1}$
$$
C_1=[B,Y].
$$
To determine the diagonal part $D_1$ we use the coefficient of $1/z$ of eq. \eqref{nwnorma}. The off-diagonal part $C_1$ does not contribute to this coefficient, so we obtain
$$
\sum_{m=1}^n{D_1}_{mm}\prod_{s\neq m}(b_s-w)=0\quad\Rightarrow\quad D_1=0.
$$

We proceed by induction. Assume that the matrices $D_1$, \dots, $D_{k-1}$, $C_1$, \dots, $C_{k-1}$ are already computed so that equations \eqref{adj}, \eqref{nwnorma} hold true modulo ${\mathcal O}(1/z^{k-1})$ and ${\mathcal O}(1/z^k)$ respectively. From the coefficient of $1/z^{k-1}$ in \eqref{adj} we have
$$
[E_a,C_{k}]= \frac{\pal C_{k-1}}{\pal x^a}+\left[ [E_a,Y], C_{k-1}\right]_{\rm off-diag} -\left[E_a,[D_{k-1},Y]\right].
$$ 
From this equation we can compute for any $i\neq a$ the $(a,i)$- and $(i,a)$-entries of the matrix $C_k$. Since $a$ is an arbitrary number between $1$ and $n$ we obtain the full off-diagonal matrix $C_k$. Equating to zero the coefficient of $1/z^k$ in \eqref{nwnorma} we obtain the diagonal matrix $D_k$
$$
{D_k}_{mm}=-\prod_{s\neq m} (b_s-b_m)^{-1} \times \mbox{coefficient of}\quad  \frac1{z^k}\quad \mbox{in}\quad 
\det\left(B+\sum_{i\leq k-1} \frac{D_i+C_i}{z^i} -b_m\cdot {\bf 1}\right).
$$
\epf

We will now prove that there exists a matrix-valued series 
\beq\label{aa}
A(\bx,z)={\bf 1}+\sum_{k\geq 1} \frac{A_k(\bx)}{z^k}
\eeq
such that
\beq\label{conj}
A^{-1}(\bx,z)M_B(\bx,z) A(\bx,z)=B
\eeq
for any diagonal matrix $B$.

\begin{lemma} For any $Y(\bx)$ satisfying eqs. \eqref{nw00} there exists a solution
\beq\label{psix}
\Psi(\bx,z)=A(\bx, z)e^{z\,\diag(x^1,\dots,\, x^n)}
\eeq
to the following system of linear differential equations
\beq\label{dpsix}
\frac{\pal \Psi}{\pal x^a}=U_{a,0}(\bx,z)\Psi(\bx,z), \quad a=1, \dots, n
\eeq
where the matrix series $A(\bx,z)$ has the form \eqref{aa}.
\end{lemma}

\pf For the coefficients $A_k=A_k(\bx)$ we obtain
$$
[E_a,A_k]=\frac{\pal A_{k-1}}{\pal x^a} +[E_a,Y]A_{k-1},\quad a=1, \dots, n.
$$
From this system we uniquely determine the off-diagonal part of the matrix $A_k$.
Using the next equation $k\mapsto k+1$ we arrive at
$$
\frac{\pal}{\pal x^a} \left(A_k\right)_{\diag}=-\left( [E_a,Y]A_{k}\right)_\diag.
$$
The off-diagonal part of $A_k$ does not contribute to the right hand side. So the matrix $A_k$ is determined uniquely up to adding a constant diagonal matrix. \epf

\begin{remark} From the proof it follows that the matrix $\Psi(\bx,z)$ is determined by eqs. \eqref{dpsix} uniquely up to a multiplication on the right by a diagonal matrix series in $1/z$
\beq\label{freex}
\Psi(\bx,z)\mapsto \Psi(\bx,z)\Delta(z), \quad \Delta(z)={\bf 1}+\sum_{k=0}^\infty\frac{\Delta^k}{z^{k+1}}, \quad \Delta^k=\diag\left(\Delta^k_1, \dots, \Delta^k_n\right)
\eeq
\end{remark}

\begin{lemma} For any diagonal matrix $B$ the solution $M_B(\bx,z)$ to the equations \eqref{adj}, \eqref{nwnorma} can be represented in the form
\beq\label{conjm}
M_B(\bx,z)=A(\bx,z) B\, A^{-1}(\bx,z)
\eeq
where the matrix $A(\bx,z)$ is defined in the previous Lemma.
\end{lemma}

\pf Since
$$
A(\bx,z) B\, A^{-1}(\bx,z)=\Psi(\bx,z) B\, \Psi^{-1}(\bx,z)
$$
the matrix \eqref{conjm} satisfies eqs. \eqref{adj}. Obviously it also satisfies \eqref{nwnorma}. Due to uniqueness of such a solution to \eqref{adj}, \eqref{nwnorma} the Lemma is proved. \epf

We are now in a position to complete the proof of Proposition \ref{prop242}. Due to uniqueness the matrix $M_B(\bx,z)$ depends linearly on $B=\diag(b_1, \dots, b_n)$. So the construction can be extended to an arbitrary diagonal matrix $B$ (see also eq. \eqref{conjm}). Put
\beq\label{ma}
M_a(\bx,z)=M_{E_a}(\bx,z), \quad a=1, \dots, n.
\eeq
These matrices clearly satisfy eqs. \eqref{adj} and \eqref{idem}. It remains to prove uniqueness.

Since the matrices $M_1(\bx,z)$, \dots, $M_n(\bx,z)$ commute pairwise due to \eqref{idem} and $M_a\to E_a$ for $z\to \infty$, we can look for their common eigenvectors in $\mathbb C^n\otimes \mathbb C\left[[ z^{-1}]\right]$. Every matrix $M_a=M_a(\bx,z)$ has only one non-zero eigenvalue; the corresponding eigenvector
$$
M_a{\bf f}_a={\bf f}_a
$$
can be normalized in such a way that $\left( {\bf f}_a\right)_b =\delta_{ab} +{\mathcal O}\left(1/z\right)$. It is determined uniquely up to multiplication
$$
{\bf f}_a\mapsto c_a(z) {\bf f}_a,  \quad c_a(z) =1+{\mathcal O}\left( \frac1{z}\right)\in \mathbb C\left[z^{-1}\right].
$$
Denote $A(\bx,z)$ the matrix whose columns are the eigenvectors ${\bf f}_1$, \dots, ${\bf f}_n$. According to the previous arguments this matrix is uniquely defined up to a multiplication on the right by $\diag (c_1(z), \dots, c_n(z))$ and satisfies
$$
M_a(\bx,z)=A(\bx,z) E_a\, A^{-1}(\bx,z), \quad a=1, \dots, n.
$$
The Proposition is proved. \epf

We will now slightly modify the setting of Proposition \ref{prop242} in order to apply it to the construction of the $n$-wave hierarchy.
Denote $\CY$ the ring of polynomials in variables $y_{ij}$, $\pal y_{ij}/\pal x^k$, $\pal^2 y_{ij}/\pal x^k \pal x^l$ etc. satisfying the constraints \eqref{nw0} along with their differential consequences. Elements of this ring will be denoted like $P(\by)$ where $P$ is a polynomial. The commuting derivations $\pal/\pal x^1$, \dots, $\pal/\pal x^n$ naturally act on this ring.

So, consider equations of the form \eqref{adj}
\beq\label{adjy}
\frac{\pal M}{\pal x^a} =\left[z E_a -[E_a, Y],M\right], \quad a=1, \dots , n
\eeq
as equations for matrices 
$$
M=M(\by,z)=M_0+\sum_{k\geq 1}\frac{M_k(\by)}{z^k}.
$$

\begin{prop} \label{propa06} There exists a unique collection of matrix series
\beq\label{may}
M_a(\by,z)=E_a+\sum_{k\geq 1} \frac{B_{a,k}(\by)}{z^k}\in Mat_n(\CY)\otimes \mathbb C\left[ z^{-1}\right], \quad a=1, \dots, n
\eeq
satisfying \eqref{adjy} and also
\beq\label{idemy}
M_a(\by,z)M_b(\by,z)=\delta_{ab}M_a(\by,z), \quad M_1(\by,z)+\dots+M_n(\by,z)={\bf 1}.
\eeq
\end{prop}

The proof essentially repeats the above arguments so it will be omitted.

\begin{remark} In practical computations of the coefficients $B_{a,k}(\by)$ instead of the normalization \eqref{nwnorma} one can alternatively use the following one:
$$
B_{a,k}(0)=0, \quad a=1, \dots, n, \quad k\geq 1.
$$
\end{remark}

Explicitly
\eqa
&&
M_a=E_a+\frac{B_{a,1}}{z}+\frac{B_{a,2}}{z^2}+\frac{B_{a,3}}{z^3}+{\mathcal O}\left(\frac1{z^4}\right)
\nn\\
&&
B_{a,1}=-[E_a,Y]
\nn\\
&&
\left(B_{a,2}\right)_{ij} =\left\{ \begin{array}{rl} -\frac{\pal y_{ij}}{\pal x^a}, & i\neq j\\
-y_{i\,a}y_{a\,i}, & j=i\neq a\\
\sum_sy_{a\, s}y_{s\,a}, & i=j=a\end{array}\right.
\nn\\
&&
\nn\\
&&
\left( B_{a,3}\right)_{ij}=\left\{ \begin{array}{rl} \frac{\pal y_{i\,a}}{\pal x^a}y_{a\, i}- y_{i\,a}\frac{\pal y_{a\,i}}{\pal x^a} , & i\neq a, ~ j\neq a\\
\\
-\frac{\pal^2 y_{a \, j}}{\pal {x^a}^2}-2 y_{a\,j} \sum_sy_{a\, s}y_{s\,a}, & i=a, ~j\neq a\\
\\
\frac{\pal^2 y_{ia}}{\pal {x^a}^2} +2 y_{i\,a}\sum_sy_{a\, s}y_{s\, a}, & i\neq a, ~ j=a\\
\\
\sum_sy_{s\,a}\frac{\pal y_{a\,s}}{\pal x^a}-\frac{\pal y_{s\,a}}{\pal x^a}y_{a\, s}, & i=j=a
\end{array}\right. .
\nn
\eeqa

Define matrix-valued polynomials
\beq\label{ualpha}
U_{a,k}(\by,z)=\left( z^{k+1} M_a(\by,z)\right)_+\in Mat_n\left(\CY\right)\otimes \mathbb C[z], \quad a=1, \dots, n, \quad k\geq -1.
\eeq
Here and below the notation $(~)_+$ will be used for the polynomial part of a Laurent series in $1/z$. 
The matrix-valued polynomials $U_{a,k}$ are exactly those that appear in the formulation \eqref{nwoper} of equations of the $n$-wave hierarchy that can be rewritten in the following form
\beq\label{nwex}
\frac{\pal Y}{\pal t^a_k}=\left( B_{a,k+2}(\by)\right)_{\rm off-diagonal}.
\eeq
For $n=2$ it coincides with the complexified nonlinear Schr\"odinger hierarchy, also known as the AKNS hierarchy. Observe that the $t^a_{-1}$-flows generate just conjugations by diagonal matrices
\beq\label{minus1}
y_{ij}\mapsto \frac{\lambda_i}{\lambda_j} \,y_{ij}, \quad i, \, j=1, \dots, n.
\eeq
Such transformations are symmetries of the $n$-wave hierarchy \eqref{nwex}.

\begin{remark} \label{rem221}The dependence of the functions $y_{ij}(\bx)$ is uniquely determined by their restriction onto any line
$$
x^i=a_i x , \quad i=1, \dots, n
$$
for arbitrary pairwise distinct constants $a_1$, \dots, $a_n$. Indeed, we reconstruct all partial derivatives in $x^1$, \dots, $x^n$ 
$$
\frac{\pal y_{ij}}{\pal x^k}=\left\{\begin{array}{cc} y_{ik}y_{kj}, & k\neq i, \, j\\
\frac{y_{ij}'}{a_i-a_j} +\sum_s \frac{a_j-a_s}{a_i-a_j}y_{is}y_{sj}, & k=i\\
\frac{y_{ij}'}{a_j-a_i}+\sum_s\frac{a_i-a_s}{a_j-a_i} y_{is} y_{sj}, & k=j
\end{array}\right.
$$
starting from the derivatives  $y_{ij}'=dy_{ij}/dx$ in $x$. So for every pair of indices $(a,k)$ the equation \eqref{nwex} can be considered as a system of $n(n-1)$ partial differential equations with one space variable $x$ and one time variable $t^a_k$.
\end{remark}

Let $Y({\bf t})$ be a solution to the $n$-wave hierarchy. Then the matrices $M_a(\by,z)$  become well-defined functions $M_a(\bt,z)$ of $\bt$. So do the matrix-valued polynomials $U_{a,p}=U_{a,p}(\bt,z)$.

\begin{prop} The matrix-valued series $M_b=M_b(\bt,z)$ satisfy
\beq\label{mrtime}
\left[ L_{a,k}, M_b\right]=0~\Leftrightarrow ~ \frac{\pal M_b(\bt,z)}{\pal t^a_k}=\left[ U_{a,k}(\bt,z),M_b(\bt,z)\right] \quad \forall~ a, \, b=1, \dots, n, \quad k\geq -1.
\eeq
\end{prop}

\pf It suffices to verify validity of eq. \eqref{mrtime} for the series $M_b=M_b(\by,z)$ and polynomials $U_{a,k}(\by,z)$. Let
$$
\tilde M_b=\frac{\pal M_b}{\pal t^a_k}-[U_{a,k},M_b].
$$
It is easy to check that this matrix-valued Laurent series satisfies
$$
\frac{\pal \tilde M_b}{\pal x^c}=\left[U_{c,0}, \tilde M_b\right]
$$
for any $c=1, \dots, n$. Let us now check that the expansion of $\tilde M_b$ contains only strictly negative powers of $z$. To this end define the matrix series
$$
V_{a,k}(\by,z)=\left(z^{k+1}M_a(\by,z)\right)_-
$$
so that
$$
z^{k+1} M_a(\by,z)=U_{a,k}(\by,z)+V_{a,k}(\by,z).
$$
Therefore
$$
\tilde M_b=\frac{\pal M_b}{\pal t^a_k}+[V_{a,k}, M_b]\in Mat_n(\CY) \otimes z^{-1} \mathbb C\left[z^{-1}\right].
$$
So, the series $\tilde M_b=\tilde M_b(\by,z)$ satisfies eqs. \eqref{adjy} and contains only strictly negative powers of $z$. Due to uniqueness it is equal to zero. \epf

\begin{lemma} The matrices $M_a(\bt,z)$ satisfy the identities \eqref{idemy}
\end{lemma}

\pf Using eq. \eqref{mrtime} we prove that 
$$
\frac{\pal}{\pal t^c_k}\left(M_a M_b-\delta_{ab}M_a\right)=0.
$$ \epf

\begin{defi} A \emph{wave function} $\Psi=\Psi({\bf t}, z)$ of the solution $Y({\bf t})$ is a solution to the infinite family of systems of linear differential equations
\beq\label{wave1}
\frac{\pal }{\pal t^a_p}\Psi=U_{a,p}({\bf t}, z) \Psi, \quad a=1, \dots, n, \quad p\geq -1
\eeq
of the form
\eqa\label{wave2}
&&
\Psi({\bf t},z)=A({\bf t}, z)e^{\phi({\bf t},z)}
\\
&&
A({\bf t},z)={\bf 1}+\frac{A^0({\bf t})}{z}+\frac{A^1({\bf t})}{z^2}+\dots,
\nn\\
&&
\phi({\bf t},z)=\sum_{k=0}^\infty\diag\left(t^1_k, \dots, t^n_k\right) z^{k+1}
\nn
\eeqa
where $A^0(\bt)$, $A^1(\bt)$ etc. are $n\times n$ matrix-valued functions of $\bt$.
\end{defi}

The wave-function is determined by a solution $Y(\bt)$ uniquely up to multiplication on the right by a constant diagonal matrix-valued series
\beq\label{wave3}
\Psi(\bt,z)\mapsto \Psi(\bt,z)\Delta(z), \quad \Delta(z)={\bf 1}+\sum_{k=0}^\infty\frac{\Delta^k}{z^{k+1}}, \quad \Delta^k=\diag\left(\Delta^k_1, \dots, \Delta^k_n\right)
\eeq

\begin{lemma} Let $\left(Y(\bt), \Psi(\bt,z)\right)$ be a solution to the equations of the hierarchy \eqref{nwex} and its wave function. Then the matrix-valued series $M_1(\bt,z)$, \dots, $M_n(\bt,z)$ can be represented in the form
\beq\label{mpsi}
M_a(\bt,z)=\Psi(\bt,z)E_a \Psi^{-1}(\bt,z), \quad a=1, \dots, n.
\eeq  
\end{lemma}

\pf As
$$
\Psi(\bt,z) E_a \Psi^{-1}(\bt,z)=A(\bt,z)E_a A^{-1}(\bt,z)=E_a+{\mathcal O}\left(\frac1{z}\right),
$$
the right hand side of \eqref{mpsi} is a series in inverse powers of $z$. It satisfies the differential equations \eqref{mrtime}. Due to uniqueness it coincides with $M_a(\bt,z)$. \epf

Introduce the following generating series for the time derivatives
\beq\label{nabla}
\nabla_a(z)=\sum_{k\geq -1} \frac1{z^{k+2}} \frac{\pal}{\pal t^a_k}, \quad a=1, \dots, n.
\eeq 

\begin{lemma} The following formula holds true
\beq\label{psit}
\nabla_a(w) \Psi(\bt,z)=\frac{M_a(\bt,w) \Psi(\bt,z)}{w-z}.
\eeq
\end{lemma}

The \emph{proof} is straightforward by using \eqref{ualpha} and \eqref{wave1}. \epf

We now proceed to the definition of tau-function. It is 
based on the following statement (cf. \cite{BDY1}).

\begin{prop} For any solution $Y(t)$ to the system \eqref{wave1} and its wave function \eqref{wave2} there exists a function $\log \tau(\bt)$ such that
\beq\label{psiz}
\frac{\pal\log\tau(\bt)}{\pal t^a_p} =-\res_{z=\infty} \tr \left( A_z(\bt, z) E_a A^{-1}(\bt,z)\right) \, z^{p+1} dz
\eeq
\end{prop}

\pf We need to prove symmetry of the second derivatives
$$
\frac{\pal^2\log\tau(\bt)}{\pal t^a_p\pal t^b_q} = \frac{\pal^2\log\tau(\bt)}{\pal t^b_q\pal t^a_p} 
$$
or, equivalently
$$
\nabla_a(z_1) \nabla_b(z_2) \log\tau(\bt)=\nabla_b(z_2) \nabla_a(z_1) \log\tau(\bt).
$$
Using eq. \eqref{psit} one can represent the generating series for the logarithmic derivatives of $\tau(\bt)$ in the following form
\beq\label{psiz1}
\nabla_a(z)\log\tau(\bt) =\tr \left( A_z(\bt, z) E_a A^{-1}(\bt,z)\right)=\tr \left( \Psi_z(\bt, z) E_a \Psi^{-1}(\bt,z)\right)-\phi_z^a(\bt,z)
\eeq
where $\phi(\bt,z)=\diag \left( \phi^1(\bt,z), \dots, \phi^n(\bt,z)\right)$ (see eq. \eqref{wave2} above).

\begin{lemma} The second order logarithmic derivatives of the tau-function \eqref{psiz}, \eqref{psiz1} can be computed from the following generating series
\beq\label{cor2}
\nabla_a(z_1) \nabla_b(z_2) \log\tau(\bt)=\frac{\tr\, M_a(\bt,z_1) M_b(\bt,z_2)-\delta_{a b}}{(z_1-z_2)^2}.
\eeq
\end{lemma}

Before we proceed to the proof let us observe that, using $M_a(z) M_b(z)=\delta_{a b}M_a(z)$ (see eq.\eqref{idemy} above) it readily follows that the numerator in \eqref{cor2} vanishes at $z_1=z_2$. Hence, due to its symmetry in $z_1$, $z_2$ it is divisible by $(z_1-z_2)^2$. Thus the right hand side is a series in inverse powers of $z_1$, $z_2$.

\pf Using the second part of eq. \eqref{psiz1} we obtain
$$
\sum_{p,\, q=0}^\infty \frac1{w^{q+2}} \frac1{z^{p+2}} \frac{\pal^2\log\tau(\bt)}{\pal t^b_q \pal t^a_p}=\nabla_b(w) \,\tr \left( \Psi_z(\bt, z) E_a \Psi^{-1}(\bt,z)\right)-\nabla_b(w) \phi^a_z(\bt, z).
$$
Obviously
$$
\nabla_b(w)\phi^a(\bt, z)=\frac{\delta_{ab}}{w-z}\quad \Rightarrow \quad \nabla_b(w)\phi^a_z(\bt, z)=\frac{\delta_{ab}}{(w-z)^2}.
$$
From \eqref{psit} obtain
$$
\nabla_b(w) \Psi_z(\bt,z)=\frac{M_b(\bt,w) \Psi_z(\bt,z)}{w-z}+\frac{M_b(\bt,w)\Psi(\bt,z)}{(w-z)^2}
$$
and
$$
\nabla_b(w) \Psi^{-1}(\bt,z)=-\frac{\Psi^{-1}(\bt,z) M_b(\bt,w)}{w-z}.
$$
Therefore
$$
\nabla_b(w) \,\tr \left( \Psi_z(\bt, z) E_a \Psi^{-1}(\bt,z)\right)=\tr\frac{M_b(\bt,w) \Psi(\bt,z)E_a \Psi^{-1}(\bt,z)}{(w-z)^2}=\tr\frac{M_b(\bt,w) M_a(\bt,z)}{(z-w)^2}.
$$
Summarizing we arrive at
$$
\sum_{p,\, q=0}^\infty \frac1{w^{q+2}} \frac1{z^{p+2}} \frac{\pal^2\log\tau(\bt)}{\pal t^b_q \pal t^a_p}=\tr\frac{M_b(\bt,w) M_a(\bt,z)}{(z-w)^2}-\frac{\delta_{ab}}{(w-z)^2}
$$
that completes the proof of Lemma and, therefore of the Proposition. \epf

\begin{remark} The definition \eqref{psiz}, \eqref{psiz1} does depend on the normalization of the wave function. A change of the normalization
$$
\Psi(\bt,z)\mapsto \Psi(\bt,z) \Delta(z), \quad \Delta(z)=\diag\left( \Delta_1(z), \dots, \Delta_n(z)\right), \quad \Delta_a(z)\in\mathbb C[[z^{-1}]], \quad a=1, \dots, n
$$
yields
$$
\tr\left( \Psi_z(\bt,z)E_a\Psi^{-1}(\bt,z)\right)\mapsto 
\tr\left( \Psi_z(\bt,z)E_a\Psi^{-1}(\bt,z)\right)+\frac{d}{dz}\log\Delta_a(z).
$$
The tau-function will change as follows
\beq\label{changetau}
\tau(\bt)\mapsto e^{\sum c_{a,p} t^a_p} \tau(\bt), \quad \frac{d}{dz}\log\Delta_a(z)=\sum\frac{c_{a,p}}{z^{p+2}}.
\eeq
\end{remark}

We see that the logarithmic derivatives of the tau-function of order two (and, therefore, of any higher order) belong to the ring ${\mathcal Y}$. In particular they do not depend on the choice of a wave function of a solution $Y(\bt)$ of the hierarchy \eqref{nwex}. Explicitly,

\beq\label{tau00}
\frac{\pal^2\log\tau}{\pal t^a_0 \pal t^b_0}=\left\{\begin{array}{cc} -y_{a b}y_{b a}, & b \neq a\\
\\
\sum_s y_{a\, s}y_{s\,a}, & b=a
\end{array}\right.
\eeq

\beq\label{tau01}
\frac{\pal^2\log\tau}{\pal t^a_0 \pal t^b_1}=\left\{\begin{array}{cc} \frac{\pal y_{a b}}{\pal x^b}y_{b a}-y_{a b} \frac{\pal y_{b a}}{\pal x^b}, & b\neq a\\
\\
\sum_s  y_{a s} \frac{\pal y_{s a}}{\pal x^s}-\frac{\pal y_{a s}}{\pal x^s}y_{s a}
, & b=a\end{array}\right.
\eeq

\beq\label{tau02}
\frac{\pal^2\log\tau}{\pal t^a_0 \pal t^b_2}=\left\{\begin{array}{ll} -y_{a b} \frac{\pal^2y_{b a}}{\pal{x^b}^2} -y_{b a}\frac{\pal^2 y_{a b}}{\pal {x^b}^2}+\frac{\pal y_{a b}}{\pal x^b}\frac{\pal y_{b a}}{\pal x^b}-3 y_{a b} y_{b a}\sum_s y_{b\, s}y_{s\, b}, & b\neq a\\
\\
-\sum_{s \neq a} \frac{\pal^2\log\tau}{\pal t^s_0 \pal t^a_2}, & b=a
\end{array}\right.
\eeq
etc.

For computation of the derivatives of order three and higher we will need the following

\begin{lemma} The following equations hold true for all $a$, $b=1, \dots, n$
\beq\label{nablam}
\nabla_a(z_1) M_b(z_2)=\frac{\left[M_a(z_1),M_b(z_2) \right]}{z_1-z_2}.
\eeq
\end{lemma}

Here and below we omit the explicit dependence on $\bt$ of the matrix-valued functions $M_a(\bt,z)$.

\pf it easily follows from \eqref{mrtime}. \epf

\begin{prop} The following equation holds true
\beq\label{cor3}
\nabla_a(z_1)\nabla_b(z_2)\nabla_c(z_3) \log\tau(\bt)=-\tr\,\frac{[M_a(z_1), M_b(z_2)] M_c(z_3)}{(z_1-z_2)(z_2-z_3)(z_3-z_1)}.
\eeq
\end{prop}

\pf it can be easily obtained by applying the operator $\nabla_c(z_3)$ at both sides of eq. \eqref{cor2} with the help of \eqref{nablam} and then using invariance of the trace of product of matrices with respect to cyclic permutations. \epf

Higher order logarithmic derivatives of the tau-function can be computed using the following

\begin{prop} For the logarithmic derivatives of order $N\geq 3$ of the tau-function of any solution $Y({\bf t})$ to the $n$-wave hierarchy \eqref{nwex} the following expression holds true
\beq\label{corN}
\sum_{k_1, \dots, k_N\geq 0} \frac{\pal^N \log \tau(\bt)}{\pal t_{k_1}^{a_1}\dots \pal t_{k_N}^{a_N}}\frac1{z_1^{k_1+2}\dots z_N^{k_N+2}}=-\frac1{N}\sum_{s\in S_N}\frac{\tr\,\left[  M_{a_{s_1}}\left( z_{s_1}\right)\dots M_{a_{s_N}}\left( z_{s_N}\right)\right]}{\left(z_{s_1}-z_{s_2}\right)\dots \left( z_{s_{N-1}} -z_{s_N}\right) \left( z_{s_N}-z_{s_1}\right)}
\eeq
\end{prop}

\pf For $N=3$ eq. \eqref{corN} coincides with \eqref{cor3}. For higher $N$ the proof is obtained by induction using \eqref{nablam}. It does not differ from the proof of a similar equation given in \cite{BDY3}, so we omit the details. \epf

  \section{Another \cite{Kac} definition of the principal tau-function}
  \setcounter{equation}{0}
  \setcounter{theorem}{0}

\begin{prop} \cite{Kac} For a given pair $\left( Y(\bt), \Psi(\bt,z)\right)$ consisting of a solution to the hierarchy \eqref{nwex} and its wave function there exists a function $\tau(\bt)$ such that
\beq\label{kac1}
\nabla_a(z) \log\tau(\bt)=\left( \frac{\pal}{\pal z} -\nabla_a(z)\right) \log \left[A(\bt, z)_{aa} \right].
\eeq
\end{prop}
It is easy to see that the diagonal entries of the matrix $A(\bt,z)$ do not depend on the variables $t^a_{-1}$. So,
according to this definition 
$
\frac{\pal\log\tau(\bt)}{\pal t^a_{-1}} \equiv 0.
$

For example,
\eqa
&&
\frac{\pal\log\tau}{\pal t^a_0}=-A^0_{aa}
\nn\\
&&
\frac{\pal\log\tau}{\pal t^a_1}=-2 A^1_{aa} +\left(A^1_{aa}\right)^2-\frac{\pal A^0_{aa}}{\pal t^a_0}.
\nn
\eeqa

\begin{defi} The function $\tau(\bt)$ will be called \emph{the principal tau-function} of the pair $\left( Y(\bt), \Psi(\bt,z)\right)$.
\end{defi}

Clearly the principal tau-function of a given pair $(Y,\Psi)$ is determined uniquely up to a nonzero constant factor.

\begin{remark} There are \cite{Kac} other tau-functions in the theory of the $n$-wave hierarchy. The principal one is selected by the following property: it is invariant with respect to diagonal conjugations
\eqa
&&
Y(\bt)\mapsto \Lambda\, Y(\bt)\,\Lambda^{-1}, \quad \Psi(\bt,z)\mapsto \Lambda\, \Psi(\bt, z)\, \Lambda^{-1}, \quad \Lambda=\diag (\lambda_1, \dots, \lambda_n),
\nn\\
&&
\tau(\bt)\mapsto \tau(\bt).
\nn
\eeqa
In other words, it does not depend on the time variables $t^a_{-1}$.
Thus
its logarithmic derivatives, starting from the second one are combinations of the functions $y_{ij}(\bt)$ and their derivatives invariant with respect to the diagonal conjugations \eqref{minus1}.
\end{remark}

{\small \it SISSA, Via Bonomea, 265, Trieste, Italy}

\end{document}